\documentstyle[12pt]{article}
\newcommand{\numero}[1]{
\addtocounter{section}{1}
\begin{center}{\bf \thesection .\ 
#1\vspace{-.1in}}\end{center}
\setcounter{subsection}{0}
\setcounter{lemma}{0}\indent}
 
\newcommand{\subnumero}[1]{
\pagebreak[1]\begin{center}{\em #1}\nopagebreak\end{center}
}
 
\newcommand{\eop}{\hfill $\Box$\vspace{.1in}}

\newtheorem{lemma}{Lemma}[section]
\newtheorem{theorem}[lemma]{Theorem}

\newtheorem{corollary}[lemma]{Corollary}

\newtheorem{conjecture}[lemma]{Conjecture}
\newtheorem{proposition}[lemma]{Proposition}

\newcommand{\cc}{{\bf C}}
\newcommand{\rr}{{\bf R}}
\newcommand{\qq}{{\bf Q}}

\newcommand{\zz}{{\bf Z}}

\newcommand{\pp}{{\bf P}}

\newcommand{\Cc}{{\cal C}}
\newcommand{\Ee}{{\cal E}}
\newcommand{\Ff}{{\cal F}}
\newcommand{\Gg}{{\cal G}}
\newcommand{\Qq}{{\cal Q}}

\newcommand{\Oo}{{\cal O}}

\newcommand{\Mm}{{\cal M}}

\newcommand{\Kk}{{\cal K}}

\newcommand{\Ll}{{\cal L}}
\newcommand{\Hh}{{\cal H}}

\newcommand{\Pp}{{\cal P}}
\newcommand{\Rr}{{\cal R}}

\newcommand{\Xx}{{\cal X}}

\newcommand{\Gm}{{\bf G}_m}

\newcommand{\gerg}{{\bf g}}

\newcommand{\gerk}{{\bf k}}

\newcommand{\germ}{{\bf m}}

\newcommand{\gerp}{{\bf p}}

\newcommand{\delbar}{\overline{\partial}}

\begin{document}

\section*{The Hodge filtration on nonabelian cohomology}

Carlos Simpson\newline
Laboratore Emile Picard\newline
UMR 5580, CNRS
Universit\'e Paul Sabatier\newline
31062 Toulouse CEDEX, France

\numero{Introduction}

Whereas usual Hodge theory
concerns mainly the usual or abelian cohomology of an algebraic variety---or
eventually the rational homotopy theory or nilpotent completion of $\pi _1$
which are in some sense obtained by extensions---nonabelian Hodge theory
concerns the cohomology of a variety with nonabelian coefficients.  Because of
the basic fact that homotopy groups in higher dimensions are abelian, and since
cohomology theories can generally be interpreted as spaces of maps into
classifying (or Eilenberg-MacLane) spaces, nonabelian cohomology occurs
essentially only in degree $1$.  There are certainly some degree $2$ aspects
which are as of yet totally untouched; and the same goes for the degree $1$
case with twisted coefficient systems. (See however \cite{kobe} for a
direction of development combining the nonabelian coefficients in degree $1$
with abelian coefficients in higher degrees).  If we leave these aside, we are
left with the case of $H^1(X,G)$ for $G$ a nonabelian group.  It is most
natural to interpret this cohomology as a groupoid, or, when $G$ is a
group-scheme, to interpret $H^1(X,G)$ as a stack.  It is the stack of flat
principal $G$-bundles on $X$.  Recall from the usual abelian
case that in order to obtain a Hodge structure, we must consider cohomology
with complex coefficients.  The analogue in the nonabelian case is that we must
take as coefficient group a group-scheme $G$ over the complex numbers (and in
fact it should be affine too).  This then is the domain of application of the
work that has been done in nonabelian Hodge theory: the study of 
properties and additional structure on the moduli
stack $\Mm (X,G):= H^1(X, G)$ which are the analogues in an appropriate sense
of the main structures or properties of abelian cohomology.

By its nature, the first nonabelian cohomology is an invariant of the
fundamental group $\pi _1(X)$.  The study of nonabelian Hodge theory may thus
be thought of as the study of fundamental groups of algebraic varieties or
compact K\"ahler manifolds.  It is important to note, specially in light of
Toledo's examples of $\pi _1(X)$ not residually finite, that the study of $\pi
_1(X)$ via its nonabelian cohomology, i.e. via the spaces of homomorphisms
$\pi _1(X)\rightarrow G$, will only ``see'' a certain part of $\pi _1(X)$ and
in particular will not at all see the intersection of subgroups of finite
index. It is an interesting question to try to understand what Hodge-theoretic
methods could say about this more mysterious part of the fundamental group.

We start in \S 2 by reviewing Corlette's nonabelian Hodge theorem
\cite{Corlette} (cf also \cite{Donaldson} and \cite{DietrichOsawa}) which is
actually a generalization of the theorem of Eells and Sampson \cite{E-S}.  This
theorem allows us to choose a prefered metric on any flat bundle.  By a Bochner
technique on K\"ahler manifolds (\cite{Siu} \cite{Jost-Yau} \cite{Corlette}) a
harmonic metric is in fact pluriharmonic, and we recover in this way the
holomorphic  data of a {\em Higgs bundle}.  The correspondence in the other
direction characterizes exactly which Higgs bundles arise in this way
\cite{Hitchin} \cite{CVHS}.  

We then mention (refering to \cite{NitsureModuli}, 
\cite{Moduli} for proofs) the
existence of moduli spaces for all of the objects in question.  The space
$M_{DR}(X,G)$ is the moduli scheme for principal holomorphic $G$-bundles with
integrable connection.  It is a coarse moduli space for the moduli stack $\Mm
_{DR}(X,G)$; it is this moduli stack which should be thought of as the
nonabelian de Rham cohomology, and the moduli space is a convenient
scheme-theoretic version.  The space $M_{Dol}(X,G)$ is the moduli scheme for
semistable principal Higgs bundles with vanishing rational Chern classes; again
it is a coarse moduli space for $\Mm _{Dol}(X,G)$, the moduli
stack which is what should be thought of as the nonabelian
Dolbeault cohomology.  The harmonic metric construction and the Bochner
technique (together with the inverse construction) give a homeomorphism
$M_{DR}(X,G) \cong M_{Dol}(X,G)$ which is $\Cc ^{\infty}$ on the smooth points.

The above work is the result of a long series of generalizations of the
original work of Narasimhan and Seshadri.  Without being exhaustive, I should
at least mention the names of Mumford, Gieseker, Maruyama, Mehta and Ramanathan
for the constructions of moduli spaces; and Donaldson, Uhlenbeck, Yau, Deligne
and Beilinson for the inverse construction to the harmonic map construction.
See the introductions and references of \cite{HBLS}, \cite{Moduli} for more
detailed historical references. 

\begin{center}
$\ast \hspace*{2cm} \ast \hspace*{2cm} \ast$
\end{center}

These constructions (or their predecessors in work of Eells and Sampson
\cite{E-S} \cite{Sampson} and Siu \cite{Siu}) are the starting point for 
much of
the work which has been done in nonabelian Hodge theory in the past several
years. As many aspects are covered by other lectures (and their corresponding
articles), and in any case many of the papers on the subject contain 
survey-like
introductions,  I will not try to survey all of these topics in the main part
of the paper but will just mention some of them here in  the introduction. 

Hitchin was interested from the beginning in the completely integrable
holomorphic hamiltonian system given by the moduli space of Higgs bundles
\cite{HitchinDuke}. This direction of research has branched off toward the
Verlinde formula, quantization and so forth. I won't try to give references as
this gets away from our principal concern of Hodge theory.

One of the principal applications of Hodge theory has always been to give
restrictions on the topological type of varieties and their subvarieties.  The
nonabelian version presents this same feature.  The  existence of all of
the various structures described above on the nonabelian cohomology $H^1(X, G)$
and various related considerations give restrictions on which groups can be
fundamental groups of compact K\"ahler manifolds, or more generally on which
homotopy types can arise.  Some of these results such as
\cite{Carlson-Toledo} \cite{Jost-Yau} \cite{Sampson} \cite{Siu} pre-date the
general Hodge-theoretic point of view given above, being based on harmonic map
considerations \`a la Eells-Sampson and Siu.  Others come directly from the
full correspondence between Higgs bundles and local systems and the subsidiary
fact that Higgs bundles invariant under the natural action of $\Gm = \cc
^{\ast}$ correspond to variations of Hodge structure \cite{HBLS}.
The restrictions on higher homotopy types are generally of two sorts: 
either one
uses information about a homotopy type to obtain additional information on the
harmonic map (such as its rank) and then concludes that such a harmonic map
cannot exist (\cite{Siu} \cite{Carlson-Toledo} \cite{CorletteRigidity}); or one
uses various notions about the cohomology of local systems to rule out higher
homotopy types (these are the restrictions coming from work of
Green-Lazarsfeld \cite{G-L}, Beauville \cite{Beauville}, Arapura
\cite{Arapura}, recently Hironaka \cite{Hironaka}, also \cite{ENS}---the idea
of using these results to get  restrictions will be explored in \S 3 
below since
many of these papers don't explicitly mention the aspect ``restrictions on
homotopy types'' which comes out of their results).

Gromov has an $L^2$-Hodge theoretic argument to rule out free (and certain
amalgamated) products of groups \cite{Gromov}.  This is particularly
interesting in relation to the theory we sketch here, because it allows one to
``see'' the whole fundamental group (for example, the amalgamated product of
two groups with no subgroups of finite index is ruled out, which would 
evidently
beimpossible to do by looking at representations into linear groups). 
Gromov and
Schoen \cite{GromovSchoen} have also developed a generalization of the 
harmonic
map theory to cover harmonic maps into negatively curved Euclidean buildings.
Coupled with a Lefschetz technique \cite{LefHolLeaves} this gives results about
fundamental groups and in particular seems to give an alternative proof of the
result about amalgamated products. This technique is generalized in
\cite{Jost-Zuo}
\cite{Katzarkov}.

With all of these restrictions on fundamental groups coming from Hodge theory,
one might well wonder if there are any interesting fundamental groups at all.
Toledo's example of a  non residually finite fundamental group
(\cite{Toledo}---there have since been several other generalizations) shows
that the family of groups which can occur is, on the contrary, quite
complicated.  We can still ask whether the part ``seen'' by nonabelian Hodge
theory with linear groups as coefficients can be nontrivial, for example, are
there nontrivial ways of obtaining positive dimensional moduli spaces (other
than easily known ways using curves and abelian varieties)?  The answer here
is affirmative too, and in fact the Higgs bundle picture is
essential for calculating what happens to obtain examples \cite{Families}. 

One of the main types of results has been the {\em factorization theorem}. 
This type of result relates the fundamental group and the geometry of $X$. 
The typical type of statement is that if $\rho : \pi _1(X)\rightarrow 
\Gamma$ is
a certain type of representation then it must {\em factor:} there is a
morphism of varieties $X\rightarrow Y$ such that $\rho$ factors through $\pi
_1(Y)$.  Perhaps the original result of this type is that of Siu \cite{Siu}
(cf also Beauville's appendix to \cite{Catanese}) stating that when $\Gamma$
is the fundamental group of a Riemann surface, then any $\rho$ must factor
through a Riemann surface $Y$.  Gromov's result of \cite{Gromov} also passes
through a similar type of statement: certain $L^2$ cohomology classes on the
universal cover of $X$ must factor through maps to a Riemann surface. One of
the first statements involving factorization through a higher dimensonal
variety is that of Zuo \cite{Zuo}; and we now have a fairly complete picture
of this type of result (cf  \cite{Jost-Zuo} \cite{Jost-Zuo2}
\cite{Katzarkov} \cite{KatzarkovPantev} \cite{Zuo} \cite{ZuoBook}): any 
nonrigid
representation to a linear group $G$ must factor through a variety $Y$ of
dimension less than or equal to the rank of $G$.  Note that the example of
\cite{Families} shows that we do not always get factorization through a
curve...but apart from this we do not know for sure if the bound rank of $G$ is
sharp. See also \cite{Mok} and various generalizations for
additional geometric information on the factorization variety $Y$.

Using the theory of harmonic maps to buildings mentioned above, one can
extend these factorization theorems to the case of representations into
linear groups over $p$-adic fields not going into a maximal compact subgroup
\cite{Jost-Zuo} \cite{Katzarkov} \cite{ZuoBook}.

These factorization techniques obviously have a certain application to the
Shafarevich conjecture \cite{Kollar}.  This has been pursued by Katzarkov,
Lasell, Napier, Ramachandran (\cite{KatzaSha} \cite{KatzaRama} 
\cite{LasellRama}
\cite{Napier} \cite{NapierRama}), see also Zuo \cite{ZuoBook}.  The main problem
is that only the part of the fundamental group seen by a given linear
representation can be treated.  They obtain a full proof of the Shafarevich
conjecture for surfaces whose fundamental group injects into a linear group.

Another recent development worth mentioning is Reznikov's proof of the Bloch
conjecture that the Chern-Simons classes of flat bundles on K\"ahler manifolds
are torsion \cite{Reznikov}.  

The last principal area of work I would like to mention, one where there is
still a fair amount to be done, is the noncompact (quasiprojective) case.  The
problem is to do the analogues of everything which we discuss in the compact
K\"ahler (or projective algebraic) case, in the case of a quasiprojective
variety.  This problem becomes much more difficult in dimension $\geq 2$. 
Aside from the dimension distinction, the problem can be divided into several
parts.  

The
first part is to obtain the analogue of Corlette's theorem.  The main
difficulty here is to get a starting point for a heat equation minimization
process, that is to say an equivariant map of finite energy.  If the
eigenvalues of the monodromy at infinity are not of norm one, this becomes
impossible and the problem is more difficult.  Modulo this difficulty, the
problem has been solved by Corlette \cite{CorletteToulouse} and
Jost and Zuo \cite{Jost-Zuo}.  

The next problem is to obtain the analogues of
the Bochner results, yielding a Higgs bundle. This is discussed to some extent
in \cite{Jost-Zuo}.  There may be a problem with Chern classes in general. 
The other aspect of this problem is that the appropriate Higgs bundle notion
must include some data at infinity, namely a {\em parabolic structure}.  
The problem of associating a parabolic Higgs bundle with nice properties, to a
harmonic bundle, is treated in Biquard \cite{Biquard} in the case when the
divisor at infinity is smooth.  

Biquard also provides  the converse
construction: given a parabolic Higgs bundle satisfying appropriate
conditions, he gets back a Yang-Mills connection and hence a representation.
This provides a relatively complete generalization in the case of smooth
divisor at infinity. What is left open for the moment is to treat the case 
where
the divisor at infinity has normal crossings.

A fourth aspect of the problem is to construct moduli spaces.  This is now
well understood for Higgs bundles and the like, due to work of Yokogawa and
Maruyama \cite{MaruYokogawa} \cite{Yokogawa}.  I think there is still a little
work (probably not too hard) left to be done on the side of {\em filtered local
systems}, which are the general representation-like objects which correspond to
parabolic Higgs bundles.

It remains to be seen how all of these aspects fit together, and then to
proceed with the generalizations of all of the further structures inherited by
the moduli spaces in the compact case (i.e. the structures we will discuss
in the present paper).

\begin{center}
$\ast \hspace*{2cm} \ast \hspace*{2cm} \ast$
\end{center}

Rather than going into further detail on all of these
applications and developments, I would like in the body of the paper to
concentrate on a fundamental aspect---the nonabelian Hodge filtration
\cite{NAHT}.   We will discuss many of the basic subjects surrounding the Hodge
filtration, such as the quaternionic structure and twistor spaces. And we give
proofs of the results announced in \cite{NAHT}, in particular the
compactification of $M_{DR}$ which is a consequence of (and practically
equivalent to) the Hodge filtration. The goal will be in the last section to
introduce an open problem, that of studying degenerations of nonabelian Hodge
structure coming from a degenerating family of varieties.

After our discussion of Corlette's theorem, its converse and the moduli spaces
in \S\S 1-2, we turn in \S 3 to a discussion of Hitchin's quaternionic 
structure
on the moduli space for representations \cite{Hitchin} \cite{Fujiki}.  
We give an application
to cohomology jump-loci retrieving the results of Green and Lazarsfeld via an
argument of Deligne and along the way see how these results give restrictions
on the higher homotopy type of non-simply connected K\"ahler manifolds.  

In \S 4 we look at Deligne's
complex analytic construction of the twistor space corresponding to the
quaternionic structure \cite{DeligneLetter}.  
He
obtains the twistor space by glueing two copies of a family $M_{Hod}$
deforming between the moduli space $M_{DR}$ of vector bundles with integrable
connection and the moduli space $M_{Dol}$ of Higgs bundles. This deformation
(parametrized by ${\bf A}^1$)
is the moduli space of {\em  vector bundles with $\lambda$-connections};
over $\lambda = 0$ a $\lambda$-connection is just a Higgs field, whereas over
$\lambda \neq 0$ a $\lambda$-connection is $\lambda$ times an integrable
connection.   

In \S 5 we explain the analogy with Rees modules which allows us to interpret
the space $M_{Hod}$ as the Hodge filtration on $M_{DR}$.  

As justification we establish in \S 6 the relationship with the Morgan-Hain
Hodge filtration on the nilpotent completion of the fundamental group
\cite{Morgan} \cite{HainKth}.

We then proceed with certain results about the Hodge filtration in the
nonabelian case, notably Griffiths transversality for
its variation in a  family,
and regularity of the Gauss-Manin connection at singular points of a family 
(\S
8). In order to do this, we first introduce in \S 7 the notion of {\em formal
groupoid} \cite{Berthelot} \cite{Illusie}.  This provides a general framework
for looking at connections, Higgs fields and so forth, and in particular allows
us to actions of these types of things on schemes rather than just vector
bundles.   

In a detour \S 9 we investigate Goldman-Millson theory \cite{GoldmanMillson}
for  the
local structure of $M_{Hod}$.  The {\em isosingularity principle} which says
that the singularities of $M_{Dol}$ are the same as those of $M_{DR}$,
generalizes to give a trivalization of $M_{Hod}$ formally along prefered
sections. This allows us to conclude, for example, that $M_{Hod}$ is flat over
${\bf A}^1$.

Then we come to a properness property of $M_{Hod}$; the limits of
$\Gm$-orbits always exist (\S 10).  This is the analogue of the classical
property of the Hodge filtration, that $F^0$ is the whole space. 

This weight property allows us to  obtain a compactification of $M_{DR}$ in \S
11, by taking
the quotient of an open set in $M_{Hod}$ by the action of $\Gm$.  
This compactification was announced without proofs in \cite{NAHT} so we take
this opportunity to provide a complete version of the argument.
(Drinfeld recently informed me that some of his students have obtained a
compactification of the moduli space of logarithmic connections on $\pp ^1$
with singularities at a finite set of points, also using the method of
$\lambda$-connections but independently of \cite{NAHT}.) 

At the end of the
section we revisit Griffiths transversality in terms of this compactification:
it says that the Gauss-Manin connection on $M_{DR}$ has poles of order $1$ at
infinity in the compactification.    This gives a picture of a compact space
$\overline{M_{DR}}(X/S,G)$ with a lift of a frame vector field on $S$ to a
vector field having a simple pole at infinity. A similar interpretation holds
for the regularity of the Gauss-Manin connection.

In the penultimate section 12 we define the {\em nonabelian Noether-Lefschetz
locus} $NL(X/S, GL(n))$.  If $X\rightarrow S$ is a family then this is
essentially the locus of $s\in S$ where $X_s$ supports an integral variation of
Hodge structure.  It is the nonabelian analogue of the classical
Noether-Lefschetz locus of Hodge cycles.  If $S$ is projective then we can see
that $NL(X/S,GL(n))$ is algebraic (as would be a consequence of the Hodge-type
conjecture that one  could
formulate, that integral variations of Hodge structure are motivic). We
conjecture that this is true even for $S$ quasiprojective, which would be a
nonabelian version of the result of \cite{CDK}.

In the last section we present an  open problem, the problem of understanding
the degeneration of all of our structures (Hodge filtration, Gauss-Manin
connection, quaternionic structure, etc.) near degenerations of a family
$X/S$.  This problem is motivated by the problem of proving that $NL(X/S)$
is algebraic when $S$ is quasiprojective---the nonabelian analogue of the
work of Cattani, Deligne and Kaplan for the classical Noether-Lefschetz
locus of Hodge cycles \cite{CDK}.

I would like to thank P. Deligne for sharing with me his ideas on how to
construct Hitchin's twistor space.  This construction provided the starting
point for everything (new) done below.

Everything is over the field $\cc$ of complex numbers.

\numero{The nonabelian Hodge theorem}

Suppose $X$ is a Riemannian manifold with basepoint $x\in X$, and suppose
$G$ is a reductive algebraic group.  A representation
$\rho : \pi _1(X,x)\rightarrow G$ corresponds to a flat principal left
$G$-bundle  over $X$ (in other words a locally constant sheaf of
principal homogeneous spaces for $G$ over $X$), or equally well to a
$\Cc ^{\infty}$ principal $G$-bundle $P$ with an integrable 
connection $D$.  We think of a connection on a principal bundle as a
$G$-invariant operator $\nabla$ from functions on $P$ to sections of
$T^{\ast}(X)|_P$ (satisfying a Leibniz rule with respect to functions pulled
back from the base). Such an operator then has a square $\nabla ^2$ from
functions on $P$ to sections of $\bigwedge ^2T^{\ast}(X)|_P$, and the
integrability condition is $\nabla ^2=0$. The flat principle bundle is the
sheaf of $\nabla$-horizontal sections of $P$.

Fix a maximal compact subgroup $K\subset G$. A {\em $K$-reduction} for a
principal bundle $P$ is a ${\cal C}^{\infty}$ principal $K$-subbundle
$P_K\subset P$ giving $P=P_K \times ^KG$.  

A flat principal $G$-bundle $(P,\nabla )$ gives rise to a flat family of
homogeneous spaces  over $X$ which we can write as $P\times ^G (G/K)$. If
$P_K\subset P$ is a $K$-reduction for $P$ then  the image of $P_K \times (eK)$
in $P_K\times ^KG$ is a smooth section of the bundle $P\times ^G (G/K)$.
This smooth section can also be thought of as a $\rho$-equivariant map
$\phi :\tilde{X}\rightarrow G/K$. We define the {\em energy} of the
$K$-reduction or of its associated equivariant map by
$$
\Ee (\phi ) := \int _X | d\phi | ^2,
$$
where the integral is taken with respect to the volume form on $X$, and
the norm of the differential $d\phi$ is measured with respect to an invariant
metric on $G/K$ (note that one has to fix a $K$-invariant metric on the
complement $\gerp $ to $\gerk \subset \gerg$ when making this discussion---if
$G$ is semisimple then we can fix the Killing form as a canonical choice).

An equivariant map or $K$-reduction $P_K$ is called {\em harmonic} if it is a
critical point of $\Ee (\phi )$.
The Euler-Lagrange equation for for a harmonic equivariant map is $d^{\ast}
d\phi =0$. This is a nonlinear equation with the Laplacian as its principal
term.

The main theorem in this subject is the following generalisation of the
theorem of Eells and Sampson \cite{E-S}:

\begin{theorem} {\rm (Corlette \cite{Corlette})}
If $\rho$ is a representation such that  $\rho (\pi
_1(X,x))$ is Zariski-dense in $G$ (or such that the Zariski closure is itself
reductive), then there exists a harmonic equivariant map $\phi$.  
\end{theorem}
{\em Proof:}
We indicate here a variant of Corlette's proof which might be useful for
people with an algebraic geometry background. Assume that the Zariski
closure $G$ is semisimple (the general reductive case may then be obtained by
using the linear theory of harmonic forms for $\cc ^{\ast}$ representations). 

Eells and
Sampson \cite{E-S} prove the existence of harmonic maps from $X$ to a compact
negatively curved manifold $M$ by a heat-equation minimization technique.  We
can start off with an equivariant map $\phi _0$ and apply the same heat
equation to obtain a family of maps $\phi _t$.  We get the same local
estimates.  In particular, for $x,y\in \tilde{X}$ the distance from $\phi _t
(x)$ to $\phi _t(y)$ is a bounded function of $t$.  Furthermore, if for any one
point $x$ we can show the ``$C^0$-estimate'' that $\phi _t(x)$ stays in a
compact subset of $G/K$ then the estimates of  \cite{E-S} will allow us to show
that the $\phi _t$ converge to  a smooth harmonic map $\phi$.   Note, for
example, that  the $C^0$-estimate is not true if the Zariski closure of $\rho
(\pi _1(X,x))$ is not reductive (in fact, the existence of an equivariant
harmonic map implies reductivity of the Zariski closure). In the original case
of \cite{E-S} the target was a compact manifold so this problem was avoided.

Now apply a little bit of geometric invariant  invariant theory to get the
$C^0$-estimate.    Choose elements $g_t\in G$ bringing us
back to the basepoint: $g_t\phi _t(x) = eK$. Fix a set of generators $\gamma
_i$ for $\pi _1(X,x)$. From \cite{E-S} we know that the distance from 
$\phi _t(x)$ to $\phi _t(\gamma _ix)$ remains bounded.  As the distance on
$G/K$ is $G$-invariant, we have
$$
d(g_t\phi _t(x), g_t\phi _t(\gamma_i x))\leq C.
$$
On the other hand the equivariance of $\phi _t$ gives $\phi _t(\gamma
x)=\rho (\gamma_i )\phi _t(x)$ so (also plugging in  $g_t\phi _t(x) = eK$) we
get 
$$
d(eK, g_t\rho (\gamma_i )g_t^{-1}eK)\leq C.
$$
Since $K$ is compact the map $G\rightarrow G/K$ is proper so the $g_t\rho
(\gamma_i )g_t^{-1}$ remain in a compact subset of $G$. 

The representation variety $R:= Hom (\pi _1(X,x), G)$ embedds in a product
$G\times \ldots \times G$ by $\rho \mapsto (\ldots, \rho (\gamma _i),
\ldots )$.  The group $G$ acts on $R$
by the {\em adjoint action} $Ad(g)(\rho )(\gamma ):= g\rho (\gamma )g^{-1}$,
and this is compatible with the above embedding via the adjoint action in each
variable of $G\times \ldots \times G$. The previous paragraph tells us
that, in our situation, $Ad(g_t)\rho $ remain in a compact subset of $R$.

The basic information from the geometric invariant theory of spaces of
representations of finitely generated groups, is that the hypothesis that 
$\rho
(\pi _1(X,x))$ is Zariski-dense in  $G$ implies that the
$Ad(G)$-orbit of $\rho$ is closed in $R$. This is well known
\cite{LubotskyMagid} but we discuss it anyway in the next two paragraphs.

If $G=GL(n,\cc )$ then $\rho$ corresponds to an irreducible $n$-dimensional
representation which we denote $V_{\rho }$.   If $V_0$ is a representation in
the closure of the orbit of $\rho$ then we have a family of representations 
$\{
V_t\}$ parametrized by $t$ in a smooth curve with 
$V_0$ being the value at a point $0$ and $V_t\cong V_{\rho}$ for $t\neq 0$.  
By
semicontinuity there is a nontrivial morphism of representations from $V_{\rho
}$ to $V_0$, but since $V_{\rho}$ is irreducible this must be an isomorphism 
and
we get that $V_0$ is in the orbit of $V_{\rho}$.

To prove this for a semisimple group
$G$ note that $G$ admits a faithful irreducible representation $V$; this 
gives a
composed representation $V_{\rho}$ of $\pi _1(X,x)$.   Since $\rho$ is
Zariski-dense, $V_{\rho}$ is irreducible. Suppose we have a family $\rho_t$ of
representations (parametrized by an affine curve) with $\rho _t \sim \rho$ for
$t\neq 0$ and $\rho _0$ different from $\rho$.  This gives a family of linear
representations $V_t$ of $\pi _1(X,x)$; as above, semicontinuity and
irreducibility of $V_{\rho}$ imply that $V_0\cong V_{\rho}$.   This then
implies that $\rho$ and $\rho_0$ are conjugate by an automorphism of $G$ where
furthermore this automorphism is a limit of inner automorphisms.  The group of
outer automorphisms being finite (hence discrete) we conclude that $\rho_0$ and
$\rho$ are conjugate by an inner automorphism, that is $\rho _0$ is in the
$Ad(G)$-orbit of $\rho$.

Now we complete the proof. The $Ad(g_t)\rho$
remain in a compact subset of the orbit $Ad(G)\rho$ because of
the fact that the orbit is closed.
Since $\rho
(\pi _1(X,x))$ is Zariski-dense in  $G$, the stabilizer of $\rho$ is just the
center of $G$, which is finite (since we have assumed that $G$ is reductive).
In particular the map $G\rightarrow Ad(G)\rho$ given by the action on $\rho$
is proper, so the $g_t$ themselves remain in a compact subset of $G$.  Finally
this implies that the $\phi (x)= g_t^{-1}eK$ remain in a compact subset of
$G/K$, which is the $C^0$-estimate we need.
\eop

{\em Remark:} Donaldson proved this theorem for rank $2$ representations
independently of Corlette \cite{Donaldson}.  It was also proved by Diederich
and Ohsawa for representations into $SL(2, \rr )$ \cite{DietrichOsawa}.
On the other hand there have since been several generalizations to the
noncompact case, for example by Corlette \cite{CorletteToulouse} and Jost and
Zuo \cite{Jost-Zuo} \cite{Jost-Zuo2}.

\subnumero{The K\"ahler case}  
Assume now that $X$ is a compact K\"ahler manifold.  
Let $\omega$ denote the K\"ahler form (of a K\"ahler metric which we choose);
let $\Lambda$ denote the adjoint of wedging with $\omega$; and let $\partial$
and $\delbar$ denote the operators coming from the complex structure.  A
holomorphic principal bundle $P$ may be considered as a $\Cc ^{\infty}$
principal bundle together with a $G$-invariant operator $\delbar$ from
functions on $P$ to sections of $\Omega ^{0,1}_X|_P$ satisfying the
appropriate Leibniz rule and $\delbar ^2 =0$. 

We say that a section, map or whatever is {\em pluriharmonic} if it is
harmonic when restricted to any locally defined smooth complex subvariety.  
This
condition is independant of the choice of metric. The classical Bochner formula
states that harmonic forms are pluriharmonic. The fundamental result about
equivariant harmonic maps on K\"ahler manifolds is just the analogue:

\begin{proposition}
If $\phi$ is a harmonic equivariant map from $\tilde{X}$ to $G/K$ then 
$\phi$ is pluriharmonic.
\end{proposition}
{\em Proof:}
See \cite{Siu} \cite{Jost-Yau} and \cite{Corlette}.
\eop

Suppose $P$ is a flat principal $G$-bundle with flat connection denoted by
$d = d' + d''$.  Suppose $P_K$ is a $K$-reduction corresponding to equivariant
harmonic map $\phi$.  We can decompose the connection into a component $d^+$
preserving $P_K$ and a component $a$ orthogonal to $P_K$.  Then decompose
according to type, $d^+ = \partial + \delbar$ and $a = \theta ' + \theta
''$.  We obtain the decompositions
$$ 
d' = \partial +
\theta '
$$ 
and
$$
d'' = \delbar + \theta ''.
$$
The components
orthogonal to $P_K$ operate on functions only via the restrictions of the
functions to the fiber, which is to say that they operate on functions via
the Lie algebra $ad (P)= P\times ^G\gerg $ of $G$-invariant vector fields on
$P$.  Thus these components are sections $\theta '$  of $ad
(P)\otimes \Omega ^{1,0}_X$ and $\theta ''$  of $ad
(P)\otimes \Omega ^{0,1}_X$. The pluriharmonic map equations translate into:
$$
\delbar ^2 = 0;
$$
$$
\delbar \theta ' + \theta ' \delbar  = 0 \;\; (\mbox{which we write}\;\;
\delbar (\theta ')=0); 
$$
and
$$
[\theta ' , \theta ' ]=0.
$$
In the last equation the form coefficients are wedged and the Lie algebra
coefficients bracketed.  

We also obtain of course the complex-conjugate equations for $\partial$ and
$\theta ''$.  These are complex conjugates in view of the hermitian or
antihermitian properties of $\partial + \delbar$ and $\theta ' + \theta ''$
respectively.

The first equation says that $(P,\delbar )$ has a structure of holomorphic
principal bundle.  This is in general different from the structure of
holomorphic  principal bundle $(P,d'' )$ which comes from the flat structure.
The second equation says that $\theta '$ corresponds to a holomorphic section
which we now denote simply by $\theta \in H^0(X, ad(P)\otimes \Omega ^1_X)$;
and the third equation says $[\theta , \theta ]= 0$.  

We define a {\em principal Higgs bundle} to be a holomorphic principal
$G$-bundle $P$ together with $\theta \in H^0(X, ad(P)\otimes \Omega ^1_X)$
such that $[\theta , \theta ]= 0$.   From the previous results, a flat 
principal
$G$-bundle $P$ with harmonic $K$-reduction $P_K$ gives a principal Higgs
bundle $(P, \theta )$. 

If $(P,\theta )$ is a principal Higgs bundle and $V$ is a
linear representation of $G$ then we obtain an associated Higgs bundle
$(E,\theta _E)$ where $E=P\times ^GV$ and $\theta _E \in H^0(X, End(E)\otimes
\Omega ^1_X)$ is the associated form associated to $\theta$. 

Recall that a Higgs bundle $(E,\theta )$ is {\em stable} if for any subsheaf
$F\subset E$ preserved by $\theta$ we have $deg (F)/r(F) < deg (E)/r(E)$ (the
notion of degree depends on choice of K\"ahler class).  Say that $E$ is {\em
polystable} if it is a direct sum of stable Higgs bundles of the same slope
(degree over rank). Say that a principal Higgs bundle $P$ is {\em polystable}
if for every representation $V$, the associated Higgs bundle is polystable. 
If the generalized first Chern classes of $P$ (corresponding to all degree one
invariant polynomials) vanish then it is enough to check this for one
faithful representation $V$ (cf \cite{HBLS} p. 86).

\begin{theorem}
\label{Corr}
Suppose $\rho : \pi _1(X)\rightarrow G$ with $G$ reductive, and suppose that
the Zariski closure of the image of $\rho$ is reductive.  Let $\Pp$ be the
associated flat bundle and $\Pp_K$ be a pluriharmonic reduction.  The structure
of principal Higgs bundle $(P,\theta )$ obtained above doesn't depend on choice
of pluriharmonic reduction $P_K$.  The principal Higgs bundle has vanishing
rational Chern classes and is polystable.  Furthermore, any polystable
principal Higgs bundle with vanishing rational Chern classes arises from a
unique representation $\rho$ in this way.   
\end{theorem}
{\em Proof:}
See \cite{Hitchin} \cite{CVHS} \cite{HBLS}.
\eop

\subnumero{Moduli spaces}

Fix a reductive complex algebraic group $G$ and a smooth projective variety 
$X$.

Let $R_{DR}(X,x,G)$ denote the moduli scheme of principal $G$-bundles with
integrable connection and frame at $x\in X$ constructed in \cite{Moduli};
similarly let $R_{Dol}(X,x,G)$ denote the moduli scheme of semistable principal
Higgs bundles with vanishing rational Chern classes with a frame at $x\in X$;
and finally let $R_{B}(X,x,G):= Hom (\pi _1(X,x), G)$ denote the space of
representations of the fundamental group in $G$. In all three cases these
schemes represent the appropriate functors.  We call these spaces the de Rham,
Dolbeault and Betti {\em representation spaces}.

The group $G$ acts on each of the representation spaces.  In all three cases,
all points are semistable for an appropriate linearized line bundle, so by
\cite{GIT} the universal categorical quotients
$$
M_{DR}(X,G):=R_{DR}(X,x,G)//G
$$
$$
M_{Dol}(X,G):=R_{Dol}(X,x,G)//G
$$
$$
M_{B}(X,G):=R_{B}(X,x,G)//G
$$
exist 
\cite{LubotskyMagid} \cite{Moduli} \cite{NitsureModuli}. They are independent
of the choice of basepoint.  The points of these quotients parametrize the
closed orbits in the representation spaces.  The closed orbit in the closure of
an orbit corresponding to a given representation is the 
{\em semisimplification}
of the representation.  Two points in a representation space map to the same
point in the moduli space if and only if their semisimplifications coincide.  

For some purposes it is useful to think about the {\em moduli stacks} instead. 
These are the stack-theoretic quotients
$$
\Mm_{DR}(X,G):=R_{DR}(X,x,G)/G
$$
$$
\Mm_{Dol}(X,G):=R_{Dol}(X,x,G)/G
$$
$$
\Mm_{B}(X,G):=R_{B}(X,x,G)/G .
$$
Properly speaking, it is  the moduli stacks which should be thought of as the
first nonabelian cohomology stacks. The moduli spaces are the
hausdorffifications or associated coarse moduli spaces for the stacks (they
universally co-represent the functors $\pi _0$ of the stacks).

We have a complex analytic isomorphism 
$R_{DR}(X,x,G)^{\rm an}\cong R_B(X,x,G)^{\rm
an}$ compatible with the action of $G$ coming from the Riemann-Hilbert
correspondence between holomorphic systems of ODE's and their monodromy
representations.  This projects to the universal categorical quotients
(\cite{Moduli} \S 5) giving  $M_{DR}(X,G)^{\rm an}\cong M_B(X,G)^{\rm
an}$ as well as to the stack quotients giving 
$\Mm _{DR}(X,G)^{\rm an}\cong \Mm _B(X,G)^{\rm
an}$.

The correspondence of Theorem \ref{Corr} gives an isomorphism between the
underlying sets of points of $M_{DR}(X,G)$ and $M_{Dol}(X,G)$, because the 
points of these spaces correspond exactly to
representations which have reductive Zariski closure (really the
corresponding de Rham or Dolbeault analogues of this notion defined using the
Tannakian formalism).  This isomorphism is a homeomorphism of underlying
topological spaces \cite{Moduli} which we thus write
$$
M_{DR}(X,G)^{\rm top} \cong M_{Dol}(X,G)^{\rm top}.
$$
Hitchin's original point of view \cite{Hitchin} was slightly different, in
that he constructed a single moduli space for all objects, and noted that it
had several different complex structures.  This amounts to the same thing if
one ignores the algebraic structures (and in fact it is difficult to say
anything concrete about the relationship between the algebraic structures). 

The
homeomorphism between the moduli spaces does not lift to a homeomorphism
between the representation spaces (cf the counterexample of 
\cite{Moduli} II, pp
38-39). I don't know what happens when we look at the stacks, so I'll give that
as a question for future research.

{\em Question:}  Does there exist a natural homeomorphism 
$\Mm_{DR}(X,G)^{\rm top} \cong \Mm_{Dol}(X,G)^{\rm top}$ inducing the previous
one on moduli spaces?

In order to attack this question one must first define the notion of the
''underlying topological space'' of a stack.

\numero{The quaternionic structure on the moduli space}

Let 
$M^{\rm sm}_{DR}(X,G)$ (resp. $M^{\rm sm}_{Dol}(X,G)$, 
$M^{\rm sm}_{B}(X,G)$) denote the
open subset of smooth points of $M^{\rm sm}_{DR}(X,G)$ (resp. 
$M^{\rm sm}_{Dol}(X,G)$,
$M^{\rm sm}_{B}(X,G)$) parametrizing Zariski-dense representations 
(the notion of
Zariski denseness makes sense for the de Rham or Dolbeault spaces using the
Tannakian point of view). Then the isomorphism  
$$
M_{DR}^{\rm sm}(X,G)^{\rm top} \cong M_{Dol}^{\rm sm}(X,G)^{\rm top} 
$$
is ${\cal C}^{\infty}$ (and even real analytic) \cite{Hitchin} \cite{Fujiki}.
Denoting by $M^{\rm sm}(X,G)$ the differentiable manifold underlying these
isomorphic spaces, we obtain two complex structures $I$ and $J$ on
$M^{\rm sm}(X,G)$ coming respectively from $M_{Dol}^{\rm sm}
(X,G)$ and $M_{DR}^{\rm sm}(X,G)$. 

The tangent space to $M^{\rm sm}_{DR}(X,G)$ (resp. $M^{\rm sm}_{Dol}(X,G)$,
$M^{\rm sm}_{B}(X,G)$) at a point
corresponding to a principal bundle $P$ is $H^1_{DR}(X, ad(P))$
(resp. $H^1_{Dol}(X, ad(P))$, $H^1_{B}(X, ad(P))$). These tangent spaces have
natural $L^2$ metrics coming from the interpretation of classes as harmonic
forms.

\begin{theorem}
\label{quater} 
{\rm (Hitchin \cite{Hitchin})} Put $K=IJ$.  Then the triple
$(I,J,K)$ is a quaternionic structure for the manifold $M^{\rm sm}(X,G)$. 
Furthermore if $g$ denotes the natural Riemannian metric on $M^{\rm sm}(X,G)$
obtained from the $L^2$ metric on the tangent space induced by the harmonic
metric (these are the same up to a constant for all structures) then $g$ is a
K\"ahler metric for each of the structures $(I,J,K)$, in other words $M^{\rm
sm}(X,G)$ becomes a hyperk\"ahler manifold.  
\end{theorem}
{\em Proof:} See \cite{Hitchin} for the result when $X$ is a curve. The 
theorem for any $X$ follows from
the corresponding theorem for a curve, using the embedding $M(X,G)\subset
M(C,G)$ for a curve $C$ which is a complete intersection of hyperplane sections
in $X$.  On the other hand Fujiki proves this for a general K\"ahler manifold
$X$ \cite{Fujiki} (where the  method of taking hyperplane sections is no longer
available).

We will
see how to calculate that $I,J,K$ form a quaternionic structure in the
proof of Theorem \ref{twistor} below.
\eop

The embedding $M(X,G)\subset M(C,G)$ obtained from a hyperplane section is
complex analytic for all structures, along any naturally defined subvariety
such as the Whitney strata of the singular locus.  Consequently the smooth
points of the underlying reduced scheme structure of these strata inherit
quaternionic (and even hyperk\"ahler) structures.

\subnumero{The twistor space}

The notion of {\em twistor space} of a quaternionic manifold $N$ is explained,
for example, in \cite{HitchinBourbaki}. Suppose $N$ is a manifold with three
integrable complex structures $(I,J,K)$ defining  a quaternionic structure on
each tangent space.  We identify $\pp ^1$ with the sphere $x^2 + y^2 + z^2 =
1$  by the stereographic projection; this gives
$$
\lambda =u+ iv \leftrightarrow (x=\frac{1-|\lambda |^2}{1+|\lambda |^2}, 
y=\frac{2u}{1+|\lambda |^2},
z=\frac{2v}{1+|\lambda |^2}).
$$

The {\em twistor space} $TW(N)$ is a
complex manifold with a  $\Cc ^{\infty}$ trivialisation 
$$
TW(N) \cong N\times \pp ^1
$$
such that the projection $TW(N)\rightarrow \pp ^1$ is holomorphic; such that
for any $n\in N$ the section $\{ n\} \times \pp ^1\subset TW(N)$ is
holomorphic; and such that for any $\lambda \in \pp ^1$ the complex 
structure on
$M\times \{ \lambda \}$ is $xI + yJ + zK$ where $(x,y,z)$ corresponds to $t$ 
via
the stereographic projection defined above (note that $(xI + yJ + zK)^2=-1$).

If we denote by $I_{u+iv}$ the complex structure corresponding to $\lambda =u+
iv\in {\bf A}^1$ then one can see using the above definitions that the formula
$$
I_{u+iv} = (1 -uK + vJ)^{-1} I (1 -uK + vJ)
$$
holds (we'll need this in the proof of Theorem \ref{twistor} below).

This definition serves to determine an almost complex structure on $TW(N)=
N\times \pp ^1$.  The almost complex structure is integrable 
\cite{HitchinBourbaki} \cite{AHS} \cite{HKLR} \cite{Salamon}; in our case we
will indicate below an explicit construction which is integrable so we can
avoid using the general integrability result.

The twistor space has various other structures, notably an antilinear
involution $\sigma$ covering the antipodal involution $\sigma _{\pp ^1}$ 
of $\pp ^1$. Define $\sigma (n,t):= (n, \sigma _{\pp ^1}(t))$.  This is
antilinear because it is antilinear in the horizontal directions (along
prefered sections) and in vertical directions because $-xI-yJ-zK$ is the
complex conjugate complex structure to $xI+yJ+zK$.

The prefered sections are by definition $\sigma$-invariant.

If $N$ is a quaternionic vector space of quaternionic rank $r$ then the twistor
space may be constructed by hand.  It the direct sum bundle $\Oo _{\pp
^1}(1)^{2r}$ over $\pp ^1$. In this case one can check that the prefered
sections are the only $\sigma$-invariant sections.

\subnumero{Application: subvarieties of $M(X,\Gm )$ defined by cohomological
conditions}

The quaternionic structure gives a nice way of looking at the results of
Green-Lazarsfeld \cite{G-L} on subvarieties defined by cohomological
conditions.  In \cite{G-L} they look at the subvarietes $\Sigma
^i_k(Pic^0)\subset Pic ^0(X)$ of line bundles $\Ll$ with $h^i(\Ll )\geq k$.
They show that these are unions of translates of subtori of $Pic ^0(X)$. 
We
look at the subvarieties $\Sigma ^i_k(M)\subset M(X,GL(n) )$
consisting of those local systems $V$ such that $h^i(X,V)\geq k$.
We have not specified whether we look at $M_B$, $M_{DR}$ or $M_{Dol}$ because
the same locus is defined in all three cases, and they correspond under the
homeomorphisms $M_B\cong M_{DR}\cong M_{Dol}$.  This is due to the fact that
if $\rho$ is a representation corresponding to vector bundle $V$ with
integrable connection and corresponding to Higgs bundle $E$ then the
interpretation of cohomology classes as harmonic forms and the K\"ahler
identities between the laplacians \cite{HBLS} gives isomorphisms
$$
H^i(X, \rho )\cong H^i_{DR}(X, V) \cong H^i_{Dol}(X, E)
$$
(see \cite{HBLS} Lemma 2.2).  Now $\Sigma ^i_k(M_{Dol})$ is a 
complex analytic subvariety of $M_{Dol}$ whereas $\Sigma ^i_k(M_{DR})$ is a
complex analytic subvariety of $M_{DR}$.  At any smooth point of the reduced
subvariety, the tangent space of $\Sigma ^i_k(M)$ is preserved by both complex
structures $I$ and $J$ of the quaternionic structure; thus at smooth points
$\Sigma ^i_k$ is a quaternionic submanifold of $M^{\rm sm}$.  This puts a big
restriction on the possibilities for $\Sigma ^i_k$ (for example, it must have
real dimension divisible by $4$ i.e. even complex dimension; the same is true
for any stratum in its Whitney stratification, for intersections of various
$\Sigma ^i_k$, etc.).

Deligne pointed out \cite{DeligneLetter} that we can use this compatibility 
with
the quaternionic structure to recover the results of \cite{G-L}. This method
(which I described as an alternative in \cite{ENS}) adds to the many various
points of view on \cite{G-L} that are now available (\cite{Arapura}
\cite{Beauville} \cite{Catanese} \cite{ENS}).  It seems worthwhile to mention
the quaternionic point of view here, since similar considerations may come 
into
play for local systems of higher rank.  

Deligne makes use of the following observation which is probably classical.

\begin{lemma}
Any locally defined smooth quaternionic subvariety of a quaternionic 
vector space is flat (i.e. a linear subspace).
\end{lemma}
{\em Proof:}  
If it were not flat, the second fundamental form would be a quaternionic
quadratic form, but an easy calculation shows that this cannot
exist. 
\eop

\begin{corollary}
\label{subtori}
If $G=\Gm$ then the $\Sigma ^i_k (M)\subset M(X,\Gm )$ are unions of
translates of subtori.
\end{corollary}
{\em Proof:}  The universal covering of $M(X,G)$ is just $H^1(X, \cc )$ and
the quaternionic structure here is linear.  Thus the previous lemma
(which is a local statement) applies to show that $\Sigma ^i_k(M)$
is flat at  smooth points of reduced
irreducible components.  A standard argument (such as in \cite{G-L}) gives the
conclusion.
\eop

This property and its generalizations (one can show that the translates of
subtori are translations by torsion points \cite{Beauville}
\cite{Catanese} \cite{Arapura} \cite{ENS}) have implications for the
topology of $X$.  One can easily fabricate examples of homotopy types such
that the corresponding jump loci are not translates of subtori.  

We give a crude version here involving additions of $2$- and $3$-cells to a
torus (one can analyze in a similar way examples made by
adding cells of any dimensions).  Put $\Gamma := \zz ^a$ with $a$
even, and put $U_1= K(\Gamma , 1)$ (which we can take as a real torus). Note
that  $$ \cc \Gamma \cong \cc [t_1, t_1^{-1},\ldots , t_a, t_a^{-1}] $$
is the Laurent polynomial ring in $a$ variables, and 
$$
M_B(U_1, \Gm )=Hom
(\Gamma , \Gm )\cong \Gm ^a.
$$ 
In fact one can canonically identify $M_B(U_1,\Gm )\cong Spec (\cc \Gamma )$
(be careful that this reasoning only works well for $G=\Gm$). Let $u\in U_1$
be the basepoint and let $U_2$ be obtained from $U_1$ by attaching $m$
$2$-spheres at $u$.  Finally let $U_3$ be obtained by attaching $\ell$
$3$-cells to $U_2$ with attaching maps $\alpha _i\in \pi _2(U_2, u)$ for
$i=1,\ldots , \ell $.  We calculate the cohomology jump loci $\Sigma
^i_k(U_3)\subset M_B(U_3,\Gm )$ for $i=2,3$. Note that $U_1\hookrightarrow U_3$
induces an isomorphism on $\pi _1$ so $M_B(U_3,\Gm ) = Spec (\cc \Gamma )$ too.
If $L$ is a rank one local system corresponding to a nontrivial representation
$\rho : \Gamma \rightarrow \Gm$ then $H^2(U_2, L)\cong L_x^m$ and 
Mayer-Vietoris gives an exact sequence 
$$
0 \rightarrow H^2(U_3, L)\rightarrow L_x^m \stackrel{A(\rho )}{\rightarrow}
L_x^{\ell } \rightarrow H^3(U_3, L)\rightarrow 0.
$$
The matrix $A(\rho )$ comes from the attaching maps: we have $\pi _2(U_2)\cong
(\cc \Gamma )^m$ so the collection $\{ \alpha _i\}$ can be considered as an
$\ell \times m$ matrix $A$ with coefficients in $\cc \Gamma$.  The matrix
$A(\rho )$ is obtained by evaluating $A$ at the algebra homomorphism $\cc
\Gamma \stackrel{\rho}{\rightarrow } \cc$.  In this case the jump loci $\Sigma
^2_k(U_3) = \Sigma ^3_{k+\ell - m}(U_3)$ are the sets of $\rho \in M_B(U_3,
\Gamma )$ where $A(\rho )$ has rank $\leq m-k$. In particular they are defined
by the ideals of $m-k$ by $m-k$ minors of $A$.  Since our choice of $\alpha
_i$ and hence of $A$ is arbitrary (except that the matrix must actually have
coefficients in $\zz \Gamma $), we can get our jump loci to be any subscheme
of $\Gm ^a$ defined by equations in $\zz \Gamma$, that is any subscheme
defined over $\zz$.

We can, for example, get the jump loci to be subschemes which are not of even
complex dimension and in any case not unions of translates of subtori---this
gives constructions of many homotopy types which cannot be the homotopy types
of complex K\"ahler manifolds. We can arrange that the jump loci do not go
through the identity representation (or even any torsion point),  in
particular this characteristic of the homotopy type will not be seen by
rational homotopy theory (we can insure that the cohomology with constant
coefficients and even the rational homotopy type are those of the real torus
e.g. an abelian variety of dimension $a/2$, or equally well those
of any complex subvariety with the same $\pi _1$).

Getting back to the result of Corollary \ref{subtori}, we can  recover the
results of Green and Lazarsfeld by looking at the Dolbeault realization.  
There
is a natural embedding $Pic ^0(X)\subset M_{Dol}(X,\Gm )$ sending a line 
bundle
$\Ll$ to the Higgs bundle $(\Ll , 0)$.  The Dolbeault cohomology of 
$(\Ll , 0)$
is just the direct sum of the $H^{i-k}(X, \Ll \otimes \Omega ^k_X)$. It is a
consequence of semicontinuity that the jump loci for a direct sum must contain
as irreducible components the jump loci for each of the factors. Thus the
irreducible components of $\Sigma ^i_k (Pic ^0)$ are among the irreducible
components of $\Sigma ^i_k(M_{Dol})$, and the conclusion of the corollary
implies the result of \cite{G-L}.

\numero{Deligne's construction of 
$TW(M^{\rm sm})$}

In \cite{DeligneLetter}
Deligne indicated a complex analytic construction of the complex
manifold $TW(M^{\rm sm})$ (the idea is based on some properties that Hitchin
established).  This is interesting because the construction given above of the
quaternionic structure on $M^{\rm sm}(X,G)$ comes from the homeomorphism
$M_{DR}(X,G) \cong M_{Dol}(X,G)$ which itself comes from the non -complex
analytic harmonic metric construction.  Of course
the trivialization $TW(M^{\rm sm}(X,G))\cong M^{\rm sm}(X,G)\times \pp ^1$
depends on the harmonic metric construction.

It turns out that Deligne's construction of $TW(M^{\rm sm}(X,G))$ is useful for
two other things that were not mentioned in \cite{DeligneLetter}: (1) it gives
an approach to defining the nonabelian analogue of the Hodge filtration on
$M_{DR}(X,G)$; and (2) it gives a way of compactifying $M_{DR}(X,G)$.  
Both of these
were announced in \cite{NAHT} but with only brief sketches of proofs.  In this
paper we will fill in the details about these two things and give some natural
extensions of these ideas.  Before doing that, we review Deligne's
construction, since it has not otherwise appeared in print (to my knowledge).

An {\em antilinear morphism $T\rightarrow T'$} between two complex analytic
spaces is a morphism of ringed spaces such that the composition $\cc
\rightarrow \Oo _{T'} \rightarrow \Oo _T$ is the complex conjugate of the
structural morphism $\cc \rightarrow \Oo _T$.  An {\em antilinear involution}
of $T$ is an antilinear morphism $\sigma : T\rightarrow T$ with $\sigma ^2=1$.

Let $\sigma _{\pp ^1}$ denote the antilinear antipodal involution of $\pp ^1$. 
If $z$ is the standard linear coordinate on ${\bf A}^1$ then $\sigma (z)=
-\overline{z}^{-1}$.  Let $\sigma _{\Gm}$ denote the restriction of $\sigma$
to $\Gm \subset \pp ^1$.

The data of a morphism of complex analytic spaces $T\rightarrow \pp ^1$
together with an antilinear involution $\sigma $ covering $\sigma _{\pp ^1}$
is equivalent to the data of a morphism $T'\rightarrow {\bf A}^1$
and an antilinear involution $\sigma '$ of $T'_{\Gm}:=T'\times _{{\bf
A}^1} \Gm$.  Given  $T'$ and $\sigma '$,
let $\overline{T}'$ denote the complex conjugate analytic space (that is the
same ringed space but with structural morphism $\cc \rightarrow \Oo
_{\overline{T}'}$ the complex conjugate of the structural morphism for $T'$).
The involution $\sigma '$ becomes a complex linear isomorphism 
$$
T'_{\Gm} \cong \overline{T}'_{\Gm},
$$
which we can use to glue $T'$ to $\overline{T}'$ to obtain $T$. By
construction $T$ comes with an antilinear involution $\sigma$ (it comes from
the tautological antilinear morphism $T'\rightarrow  \overline{T}'$ and its
inverse). One can see that $\pp ^1$ is obtained from ${\bf A}^1$ and the
involution $\sigma _{\Gm}$ by the same construction, so we obtain our map
$T\rightarrow \pp ^1$ compatible with involutions.

To give a holomorphic $\sigma$-invariant section $\eta :\pp ^1\rightarrow T$ 
it
suffices to give a holomorphic section  $\eta ' : {\bf A}^1 \rightarrow T'$
such that $\eta ' |_{\Gm}$ is $\sigma '$-invariant.

Hitchin noticed that the twistor space $TW(M^{\rm sm})$  comes equipped with
an action of $\Gm$ identifying the fibers over all different $\lambda \in \Gm
\subset \pp ^1$  \cite{Hitchin} (note however that the twistor space of a
general hyperk\"ahler manifold doesn't
 come equipped with such an action).  Deligne's idea is to use this and the
remark of the previous paragraph to obtain a direct construction of 
$TW(M^{\rm sm})$.

For simplicity we treat the case
$G=GL(n,\cc )$ but the case of a general reductive group can be treated
directly by working with principal bundles, or indirectly using the Tannakian
formalism such as in \cite{Moduli} (note that for the constructions of moduli
spaces the indirect Tannakian method is the only one I know of).

Deligne makes the following definition.
Suppose $\lambda : S\rightarrow {\bf A}^1$ is a morphism.  A {\em
$\lambda$-connection} on a vector bundle $E$ over $X\times S$ consists of an
operator 
$$
\nabla : E\rightarrow E\otimes _{\Oo} \Omega ^1_{X\times S/S}
$$
such that $\nabla (ae) = \lambda e\otimes d(a) + a\nabla (e)$ (Leibniz
rule multiplied by $\lambda$) and such that $\nabla ^2 = 0$ as defined in the
usual way (integrability).   

Note that if $\lambda = 1$ then this is the same as the usual notion of
connection, whereas if $\lambda = 0$ then this is the same as the notion of
Higgs field making $(E,\nabla )$ into a Higgs bundle \cite{Hitchin}
\cite{HBLS}.

\begin{proposition}
\label{ConstructMhod}
Fix $x\in X$.
The functor which to $\lambda : S\rightarrow {\bf A}^1$ associates the set of 
triples $(E,\nabla , \beta )$ where $E$ is a vector bundle on $X\times S$,
$\nabla$ is a $\lambda$-connection on $E$ (such that the resulting Higgs 
bundles over
$\lambda=0$ are semistable with vanishing rational Chern classes), and 
$\beta : E|_{\{ x\} \times S}
\cong \Oo _S^n$ is a frame, is
representable by a scheme $R_{Hod}(X,x,GL(n)) \rightarrow {\bf A}^1$.  The 
group
$GL(n)$ acts on $R_{Hod}(X,x,GL(n))$ by change of frame and  all points are
semistable for this action (with respect to the an appropriate linearized
bundle).  The geometric-invariant theory quotient $M_{Hod}(X,GL(n))\rightarrow
{\bf A}^1$ is a universal categorical quotient. In particular the fibers of
$M_{Hod}$ over $\lambda = 0$ and $\lambda = 1$ are $M_{Dol}(X,GL(n))$ and 
$M_{DR}(X,GL(n))$
respectively.    
\end{proposition} 
{\em Proof:}
This follows by applying the results of \cite{Moduli}
to the ring $\Lambda ^R$ defined in \cite{Moduli} p. 87.
\eop

{\em Remark:}
We can define the notion of $\lambda$-connection on a principal
bundle, and obtain the corresponding statement for principal $G$-bundles. 
We obtain schemes $R_{Hod}(X,x,G)$ and $M_{Hod} (X,G)$.
The construction is done  by applying
the Tannakian considerations of \cite{Moduli} \S 9.

Concerning the terminology $R_{Hod}$ and $M_{Hod}$: 
this reflects the fact that,
as we shall see below, these spaces incarnate the Hodge filtrations on
$R_{DR}$ and $M_{DR}$.

Let $\Mm _{Hod}(X,GL(n))$ (or $M_{Hod}(X,G)$) denote the stack-theoretic
quotient of $R_{Hod}(X,x,GL(n))$ by $GL(n)$ (or $R_{Hod}(X,x,G)$ by $G$).

The group $\Gm$ acts on the functor $\{ (E,\nabla , \beta )\}$ over its action
on ${\bf A}^1$: if $t\in \Gm (S)$ and $(E,\nabla , \beta )$ is a
$\lambda$-connection then $(E, t\nabla , \beta )$ is a $t\lambda$ connection. 
Since $R_{Hod}(X,x,GL(n))$ represents the functor, we get an action of $\Gm$ on
$R_{Hod}(X,x,GL(n))$ covering its action on ${\bf A}^1$.  Since $M_{Hod}
(X,GL(n))$ is a universal
categorical quotient, this descends to an action on $M_{Hod}(X,GL(n))$.  
This action
serves to identify the fibers over any $\lambda , \lambda ' \neq 0$ in ${\bf
A}^1$---they are all isomorphic to $M_{DR}(X,GL(n))$.

The space $M_{Hod}(X,GL(n))$ will play the role of the space $T'$ in 
constructing the
twistor space $T$.  According to our general discussion, in order to obtain
$T$ by glueing, it suffices to have an antilinear involution $\sigma '$ of
$M_{Hod}(X,GL(n))|_{\Gm}$. As we have seen above, the action of $\Gm$ gives an
isomorphism 
$$
M_{Hod}(X,GL(n))|_{\Gm} \cong M_{DR}(X,GL(n)) \times \Gm .
$$
On the other hand we have an antilinear involution $\tau $ of $M_B(X,GL(n))$ 
obtained
by setting $\tau (\rho )$ equal to the dual of the complex conjugate
representation (where complex conjugation is taken with respect to the real
structure $GL(n,\rr )$; the dual of the complex conjugate is also the complex
conjugate with respect to the compact real form). To be totally explicit, for
$\gamma \in \pi _1(X,x)$ we set $\tau (\rho )(\gamma ):= \;\; ^t\overline{\rho
(\gamma )}^{-1}$. The complex analytic isomorphism $M_B(X,GL(n))^{\rm an}\cong
M_{DR}(X,GL(n))^{\rm an}$ given by the Riemann-Hilbert correspondence 
allows us to
interpret $\tau$ as an antilinear involution of $M_{DR}(X,GL(n))$.  Finally  
we define
the involution $\sigma '$ of $M_{DR}(X,GL(n)) \times \Gm $ by the formula 
$\sigma '
(u,\lambda )= (\tau (u), -\overline{\lambda}^{-1})$. Using $\sigma '$ and
$T'=M_{Hod}(X,GL(n))$ in the general recipe given above, we obtain a space $T$ 
which we
denote $M_{Del}(X,GL(n))\rightarrow \pp ^1$.

The complex conjugate scheme $\overline{T'}= \overline{M_{Hod}(X,GL(n))}$
which appears above can be identified with $M_{Hod}(\overline{X}, GL(n))$.

{\em Exercise:} write down the glueing isomorphism between  
$M_{Hod}(X, GL(n))$ and $M_{Hod}(\overline{X}, GL(n))$ over $\Gm
\subset {\bf A}^1$.  Note that it will be analytic but not algebraic
(depending on the Riemann-Hilbert correspondence).

We note rapidly some properties of $M_{Del}(X,GL(n))$ which are immediate 
consequences
of the construction.  The fiber of $M_{Del}(X,GL(n))\rightarrow \pp ^1$ 
over a point
$\lambda \in \pp ^1$ (which is denoted using the coordinate system of the
first embedding of ${\bf A}^1$ which corresponds to the part concerning $X$)
is equal to $M_{Dol}(X,GL(n))^{\rm an}$ if $\lambda =0$; the fiber is 
isomorphic to
$$
M_{DR}(X,GL(n))^{\rm an}\cong M_{B}(X,GL(n))^{\rm an}\cong
M_{DR}(\overline{X} ,GL(n)))^{\rm an} $$
if $\lambda \neq 0,\infty $; and the fiber is equal to
$M_{Dol}(\overline{X},GL(n))^{\rm an}$ if $\lambda = \infty$.
There is an analytic action of $\Gm$ covering the standard action on $\pp ^1$
(this action is constructed by glueing the natural action over the first open
set $M_{Hod}(X,GL(n))^{\rm an}$ with the composition with $i$ of the natural
action on the second open set $M_{Hod}(\overline{X},GL(n))^{\rm an}$).
The antilinear involution $\sigma $ of $M_{Del}(X,GL(n))$ comes from the 
first version
of the construction discussed above.

{\em Remark:}  Let $R_{Hod}(X,x,GL(n))$ denote
the representation space constructed with reference to the basepoint $x\in X$. 
We can construct,
exactly as above,  an involution $\sigma$ which can be thought of as an
isomorphism between the inverse images of $\Gm \subset {\bf A}^1$ in
$R_{Hod}(X,x,GL(n))$ and  $R_{Hod}(\overline{X},
\overline{x},GL(n))$.
We obtain $R_{Del}(X,x,GL(n))\rightarrow \pp ^1$ by glueing 
$R_{Hod}(X,x,GL(n))$ to 
$R_{Hod}(\overline{X},
\overline{x},GL(n))$ using this isomorphism. The group $GL(n,\cc )$ acts
analytically and on each open subset the associated moduli space
$M_{Hod}(X,GL(n))$ is a universal categorical quotient in the analytic category
(\cite{Moduli} \S 5). The glueing (along invariant open sets which are
pullbacks of open sets in the quotients) preserves this property, so 
$$
R_{Del}(X,x,GL(n))\rightarrow M_{Del} (X,GL(n)) 
$$
is a universal categorical quotient by the action of $GL(n , \cc )$ in the
analytical category.  There is again an action of $\Gm$ on 
$R_{Del}(X,x,GL(n))$ and the fibers are again respectively 
$R_{Dol}(X,x,GL(n))^{\rm an}$, $R_B(X,x,GL(n))^{\rm an}$ and 
$R_{Dol}(\overline{X},
\overline{x},GL(n))$ over $\lambda = 0$, $\lambda \neq 0,\infty$, and 
$\lambda =
\infty$ in $\pp ^1$.

Denote by $\Mm _{Del}(X,GL(n))$ the stack-theoretic quotient of 
$R_{Del}(X,x,GL(n))$ by $GL(n)$; it is an analytic stack with a morphism to
$\pp ^1$.

We now show how a {\em harmonic bundle} defines a section  
$\pp ^1\rightarrow M_{Del}(X,GL(n))$ which we refer to as a {\em prefered
section}. As mentioned before, it suffices to obtain a $\sigma '$-invariant
section ${\bf A}^1\rightarrow M_{Hod}(X,GL(n))$.

Suppose $P$ is a flat principal $GL(n)$-bundle.  Choose a pluriharmonic
$K$-reduction $P_K$ and consider the decomposition defined previously
$$
d' = \partial + \theta ',
$$
$$
d'' = \delbar + \theta ''.
$$
For $\lambda \in {\bf A}^1$ we define a holomorphic structure
$$
\delbar _{\lambda} := \delbar + \lambda \theta '',
$$
and an operator
$$
\nabla _{\lambda} := \lambda \partial + \theta ' .
$$
We claim that $\delbar _{\lambda}$ is an integrable holomorphic structure and
$\nabla _{\lambda}$ an integrable holomorphic $\lambda$-connection on
$(P,\delbar _{\lambda})$. The equations $\delbar ^2=0$, $[\theta '' , 
\theta ''
]=0$  and $(d'')^2=0$ imply that $\delbar (\theta '' )=0$ and hence $\delbar
_{\lambda}^2=0$. Similarly, $\delbar (\theta ')=0$ and $\partial (\theta
'')=0$ and furthermore we get 
$$
\delbar \partial + \partial \delbar + \theta ' \theta '' + \theta '' \theta
' = 0
$$
which gives $[\delbar _{\lambda} ,\nabla _{\lambda}]=0$; finally by the
same argument as previously $\partial (\theta ')=0$ so $\nabla _{\lambda
}^2=0$.  This gives the claim. Note that at $\lambda = 0$ we recover the
Higgs bundle structure $(\delbar , \theta ')$ which we know to be
polystable with vanishing Chern classes.  This construction thus gives a
section ${\bf A}^1\rightarrow M_{Hod}(X,GL(n))$.  It is holomorphic in 
$\lambda$
(since $\lambda$ appears linearly in the equations).  

We have to check that our section is $\sigma '$-invariant over $\Gm \subset
{\bf A}^1$.  This is a bit technical so feel free to skip it!  A point of
$M_{Hod}(X,GL(n))$ can be represented as a quadruple $(E,\delta ' , \delta '',
\lambda )$ where $E$ is a $\Cc ^{\infty}$ bundle, $\lambda \in \cc$, $\delta $
is an operator satisfying Leibniz' rule for $\lambda \partial$, and
$\delta ''$ is an operator satisfying Leibniz' rule for $\delbar$, such
that $(\delta ' )^2=0$, $(\delta '' )^2=0$, and $\delta ' \delta '' +
\delta '' \delta '=0$.  If $\lambda \neq 0$ this corresponds to a flat
bundle $(E, \lambda ^{-1}\delta ' + \delta '' )$.  The dual complex conjugate
flat bundle (corresponding to the dual of the complex conjugate representation
on $X$) is $(\overline{E}^{\ast},\overline{\delta ''}^{\ast} +\overline
{\lambda}^{-1}
\overline{\delta '}^{\ast})$ (the superscript $\; ^{\ast}$ on the operators
means the induced operators on the dual). If we take the point obtained by
multiplying this complex conjugate flat bundle by $-\overline{\lambda }^{-1}$,
we obtain the point 
$$ \sigma (E,\delta ', \delta '', \lambda ) = (\overline{E}^{\ast},
-\overline{\lambda}^{-1} \overline{\delta ''} ^{\ast}, 
\overline{\lambda}^{-1} \overline{\delta '}^{\ast}, -\overline{\lambda}^{-1}).
$$
We have to check that this operation preserves our preserved section, which is
the collection of points of the form $(E, \lambda \partial + \theta ',
\delbar + \lambda \theta '', \lambda )$.   We have an isomorphism
$\overline{E}^{\ast} \cong E$ given by the harmonic metric, and via
this isomorphism the dual complex conjugation operation has the
following effect on operators:
$$
\partial \leftrightarrow \delbar ,
$$
$$
\theta ' \leftrightarrow - \theta ' 
$$
(this is from the definition of $\delbar + \partial $ and $\theta '+\theta ''$
as the components parallel to and perpendicular to the unitary structure). In
view of these formulae, when we apply the operation $\sigma$ to such a point we
get  
$$ 
\sigma (E, \lambda \partial + \theta ', \delbar + \lambda \theta '',
\lambda ) = (\overline{E}^{\ast}\cong E,
 -\overline{\lambda}^{-1}\partial + \theta ', \delbar
-\overline{\lambda}^{-1} \theta '',  -\overline{\lambda}^{-1} )
$$
which is indeed a point on our prefered section.

It is clear from the definition that through any point of $M_{Del}(X,GL(n))$ 
passes
exactly one prefered section.

The set of prefered sections gives a set-theoretic trivialization
$$
M_{Del}(X,GL(n))\cong M_B(X,GL(n))\times \pp ^1.
$$
This trivialisation is in fact a homeomorphism, as can be seen by using the
techniques of \cite{Moduli}  which are used in proving that
$M_{DR}(X,GL(n))  ^{\rm top} \cong M_{Dol}(X,GL(n))^{\rm top}$. This is also
verified in \cite{Fujiki}.

Let  $M_{Del}^{\rm sm}(X,GL(n))$ denote the open subset of $M_{Del}(X,GL(n))$ 
where the
projection to $\pp ^1$ is smooth.  By the etale local triviality of $M_{Hod}$
(explained in \S 9 below) a point in $M_{Del}(X,GL(n))$ lies in 
$M_{Del}^{\rm sm}(X,GL(n))$ if and only if it is a smooth point of the fiber
$M_{Dol}(X,GL(n))$, $M_{DR} (X,GL(n))$ or $M_{Dol}(\overline{X}, GL(n))$.  

The trivialisation via
prefered sections gives 
$$
M_{Del}^{\rm sm}(X,GL(n))^{\rm top} \cong 
M^{\rm sm}(X,GL(n))^{\rm top} \times \pp ^1.
$$
This is in fact a ${\cal C}^{\infty}$ isomorphism, as  
follows from the construction of the prefered sections and the fact that
the harmonic maps or metrics  vary smoothly with parameters (since they are
solutions of the appropriate kind of nonlinear elliptic equation).

\begin{theorem} 
\label{twistor}
{\rm (Deligne)}
The space $M_{Del}^{\rm sm}(X,GL(n))$ with all of its structures is 
analytically
isomorphic to the twistor space $TW(M^{\rm sm})$; via this isomorphism,
the prefered section trivialisations of $M_{Del}^{\rm sm}(X,GL(n))$  and
$TW(M^{\rm sm})$ coincide.
\end{theorem}
{\em Proof:}
This is actually a consequence of the properties obtained by Hitchin for his
twistor space in \cite{Hitchin}.  For intrepid readers, we indicate a
self-contained calculation---partly because this also serves to show that
$(I,J,K)$ defined a quaternionic structure in the first place. 
Both the twistor
space and $M_{Del}^{\rm sm}(X,GL(n))$ are $\Cc ^{\infty}$ isomorphic to the
product $M^{\rm sm}\times \pp ^1$. This gives the isomorphism between the two.
Furthermore we know in both cases that the horizontal  sections $\{ x\} \times
\pp ^1$ are holomorphic, so the isomorphism is analytic in the horizontal
direction. We have to check that this isomorphism is compatible with the
complex structures in the vertical direction.   Choose a tangent direction to
$M^{\rm sm}$ which we will look at first in the Dolbeault realization.  The
tangent direction can be thought of as a change of operator $\delbar + \theta '
\mapsto \delbar + \theta ' + \alpha $ where $\alpha $ is an endomorphism-valued
form representing the cohomology class of the tangent vector.  We may
(by gauging back if necessary) assume that the associated harmonic metric
remains fixed; the infinitesimal change $\alpha$ then induces a change of
operator $\partial + \theta '' \mapsto \partial + \theta '' + \beta $.  Write
$\alpha = \alpha ' + \alpha ''$ and $\beta = \beta ' + \beta ''$ according to
type.  In
terms of the isomorphism $E\cong \overline{E}^{\ast}$ the condition that the
fixed metric still relates our new operators is 
$$
\beta '= \overline{\alpha ''}^{\ast}
$$
$$
\beta '' = -\overline{\alpha '}^{\ast}.
$$
One can see that if $\alpha$ is harmonic then the form $\beta$
defined by these formulas is also harmonic.   If we denote by $B(\alpha )$ the
form $\beta$ defined by these formulas then $B$ becomes an endomorphism of the
space of harmonic forms.  It is antilinear (that is $Bi=-iB$), and $B^2=-1$. 
Thus $B$ is another complex structure which forms part of a quaternionic triple
with $i$.  

The complex structure $I$ on $M_{Dol}(X,GL(n))$ corresponds to 
multiplication of
$\alpha$ by $i$ (because $\alpha$ is the representative of our tangent vector
in the Dolbeault realization).  The complex structure on $M_{DR}(X,GL(n))$ is
the operator on $\alpha$ which causes $\alpha + \beta = \alpha + B(\alpha )$ 
to
be multiplied by $i$. Thus we have the formula 
$$ 
I(1+ B)\alpha = (1+B)J\alpha .
$$
From whence $J=IB$.  This now shows that the pair $(I,J)$
form a part of a quaternionic triple, for which $B=-K$ (Theorem
\ref{quater}).    

For $\lambda \in {\bf A}^1$ the change of associated
$\lambda$-connection is 
$$
\lambda \partial + \theta ' +
\delbar + \lambda \theta ''\mapsto \lambda \partial + \theta ' +
\delbar + \lambda \theta ''+ \lambda \beta ' + \alpha ' + \alpha '' + \lambda
\beta '' .
$$ 
Thus if $I_{\lambda}$ denotes the complex structure on the fiber
of $M_{Del}(X,GL(n))$ over $\lambda \in {\bf A}^1$ then we get the formula
$$
I(1+ \lambda B) = (1+ \lambda B)I_{\lambda}.
$$
One has to be careful about what $\lambda B$ means: if $\lambda = u+ iv$ then 
$\lambda B = uB + vIB$.  We obtain (replacing $B$ by $-K$):
$$
I_{u+iv} = (1 -uK + vJ)^{-1} I (1 -uK + vJ).
$$
This coincides with the formula given in \S 3.
\eop

{\bf Question:} Are the prefered sections the only sections which are
preserved by the involution $\sigma$?

This is certainly locally true, since the normal bundle to a prefered
section is a direct sum of $\Oo _{\pp ^1}(1)$.  In fact, locally the morphism
from the space of all sections to the product of any two distinct fibers is
an isomorphism. If we take two antipodal fibers then $\sigma$ gives an
antilinear involution of the product of the two fibers, and the prefered
sections correspond to the fixed points.

Hitchin's discussion in \cite{HitchinBourbaki} (Theorem 1) is actually a bit
unclear on this point: as written the converse in Theorem 1 would imply that
the answer is yes in general, but one can easily imagine that he meant only to
look at the real sections in the given family of sections. A glance at
\cite{HKLR} didn't resolve the problem, so I think that the answer to the
above question is not known.

An affirmative answer would mean that $(M_{Del}(X,GL(n)), \sigma )$ determines
the twistor space structure and in particular
the isomorphism $M_{DR} \cong M_{Dol}$, an interesting point since the
construction of $(M_{Del}(X,GL(n)), \sigma )$ is entirely complex analytic,
so we could bypass the nonlinear elliptic theory necessary to define the
harmonic metrics---conceptually speaking at least.

The answer to this question is `yes' for the twistor space of a
quaternionic vector space.  As a consequence we obtain this property for the
moduli space of rank one representations:  

\begin{theorem}
Suppose $G=\Gm$.  Then the prefered sections are the only $\sigma$-invariant
sections of $TW(M^{\rm sm})\rightarrow \pp ^1$.
\end{theorem}
{\em Proof:}
The moduli space is a quotient $M = H^1(X, \cc )/H^1(X, \zz )$ (as can be
seen by a flat version of the exponential exact sequence).  The quaternionic
structure is the quotient by the lattice of a linear quaternionic structure on
$H^1(X, \cc )$. Thus the twistor space is the quotient
$$
TW(M)= TW(H^1(X, \cc )) /H^1(X, \zz ).
$$
Since $\pp ^1$ is simply connected, the sections from $\pp ^1$ to $TW(M)$ are
just projections of sections from $\pp ^1$ to $TW(H^1(X, \cc )$.  The
involution $\sigma$ acts compatibly on everything.  From the theory of the
twistor space for quaternionic vector spaces (which is just a bundle which is
a direct sum of $\Oo _{\pp ^1}(1)$) we see that through any point of
$TW(H^1(X, \cc ))$ there is a unique $\sigma$-invariant section; this gives
the same result on $TW(M)$ which implies the theorem.
\eop

\numero{The Hodge filtration}

In \cite{Hitchin} Hitchin introduced an $S^1$ action on the moduli space of
representations.  This was taken up again in \cite{HBLS} as a $\cc
^{\ast}$-action. This action is defined via the isomorphism $M_B^{\rm top}
\cong M_{Dol}^{\rm top}$: $t\in \cc ^{\ast}$ sends the Higgs bundle $(E,\theta
)$ to $(E, t\theta )$.

The $\cc ^{\ast}$ or $S^1$ actions are the analogue in
nonabelian Hodge theory of the Hodge decomposition of cohomology coming from
harmonic forms.   In the usual case, the Hodge decomposition does not vary
holomorphically with parameters, because it includes complex conjugate
information.  Similarly, if the variety is defined over a small field, there is
no particular reason for the Hodge decomposition to be defined over a small
field.  In order to obtain something which comes from algebraic geometry and
thus has the properties of holomorphic variation, and compatibility with fields
of definition, one looks at the {\em Hodge filtration} of the algebraic de Rham
cohomology.  We will define and investigate the analogue for nonabelian
cohomology. 

Begin with the following  observation. Suppose $V$ is a vector space with
complete decreasing filtration $F^{\cdot}$ (complete means that the
filtration starts with $V$ and ends with $\{ 0\}$).  Define a locally free
sheaf $\xi (V,F)$ over ${\bf A}^1$ with action of $\Gm$  as
follows.  Let $j: \Gm \rightarrow {\bf A}^1$ denote the inclusion.  Then $\xi
(V,F)$ is the subsheaf of $j_{\ast}(V\otimes \Oo _{\Gm})$ generated by the
sections of the form $z^{-p}v_p$ for $v_p \in F^pV$ (where $z$ denotes the
coordinate on ${\bf A}^1$). Conversely if $W$ is a locally free sheaf on ${\bf
A}^1$ with action of $\Gm$ then we obtain a decreasing filtration $F$ on the
fiber $W_1$ of $W$ over $1\in {\bf A}^1$ by looking at orders of poles of
$\Gm$-invariant sections. These constructions are inverses.

The locally free sheaf $\xi (V,F)$ is the tilde of the {\em Rees module} of
$(V,F)$.

If $(V,F)$ is a filtered vector space then the fiber $\xi (V,F) _0$ over $0\in
{\bf A}^1$ is naturally identified with the associated-graded $\bigoplus
F^p/F^{p+1}$.

Let $\Sigma$ be a sheaf of sets  on the big
etale site $\Xx$.  We define a {\em filtration $\Ff$ of $\Sigma$} to be a sheaf
of sets with morphism $\Sigma  _{\Ff}\rightarrow {\bf
A}^1$  together with action of $\Gm$ (here an action means
a morphism $\Sigma _{\Ff} \times \Gm \rightarrow \Sigma _{\Ff}$ satisfying 
the usual
axioms) and an isomorphism $\Sigma _{\Ff} \times _{{\bf A}^1} \{ 1\} \cong
\Sigma$.  Note that $\Sigma _{\Ff}$ may be interpreted as a sheaf on $\Xx
/{\bf A}^1$, and using this interpretation we can make a similar definition
for sheaves of objects of any appropriate category.  We obtain a similar
definition for stacks (or homotopy-sheaves of spaces, or even $n$-stacks or
$\infty$-stacks once those are defined).  

Normally we will be interested in the case where $\Sigma$ is represented by a
scheme or eventually an algebraic stack, in this case we expect $\Sigma
_{\Ff}$ to be a scheme or at least an algebraic stack.

{\em Caution:}  As we will see in one of our main examples in the section on
formal categories, the notion of filtration of a sheaf of sets in the context
of stacks is different from the notion of filtration in the context of sets,
in other words we might have $\Sigma _{\Ff}$ a stack whereas $\Sigma$ is a set.

Now, getting back to our main discussion, in terms of this definition the
space $M_{Hod}\rightarrow {\bf A}^1$ with action of $\Gm$ is a filtration on
$M_{DR}$.  We call this the {\em
Hodge filtration} on $M_{DR}$. Similarly $R_{Hod}$ is the Hodge filtration on
$R_{DR}$. And most properly speaking, it is $\Mm _{Hod}$ which
provides the Hodge filtration on the nonabelian cohomology stack $\Mm _{DR}$.

In the next section we will see how this filtration is compatible with the
usual Hodge filtration on the nilpotent completion of the fundamental group.

The idea of interpreting the Hodge filtration in this way is very closely
related to the interpretation of Deninger \cite{Deninger}. Essentially he
looks at a derivation expressing the infinitesimal action instead of the full
action of $\Gm$. In turn he refers to Fontaine \cite{Fontaine} (1979) for a
reworking of Hodge theory from this point of view (which is what led Fontaine
to all of his rings such as $B_+^{\rm cris}$\ldots I guess\ldots ).

A word about purity.  If $(V,F,\overline{F})$ is a vector space with two
filtrations (which can be complex conjugates with respect to a real structure,
for example) then $\xi (V, F, \overline{F})$ is a vector bundle over $\pp
^1$ obtained by glueing $\xi (V,F)$ to $\xi (V,\overline{F})$ much as in \S 4. 
The two filtrations define a Hodge decomposition pure of weight $w$ if and only
if the vector bundle $\xi (V, F, \overline{F})$ is a direct sum of copies of
$\Oo _{\pp ^1}(w)$.   
The construction $M_{Del}$
is in effect the nonabelian
analogue of the construction $\xi (V,F,\overline{F})$ where $F$ is replaced by
the ``filtration'' $M_{Hod}$. The fact that this construction gives the
twistor space for a quaternionic structure is equivalent to the statement that
the normal bundle along any prefered section is a direct sum of $\Oo _{\pp
^1}(1)$ (cf \cite{HitchinBourbaki}, \cite{HKLR}).  This can be interpreted as
purity of weight one. I don't know how far one can go toward making this
analogy more precise than it is.

\numero{The nilpotent completion of $\pi _1$ and representations near the
identity}

We will justify our definition of $M_{Hod}$ as the Hodge filtration by making
the connection with the usual Hodge filtration on the nilpotent completion of
the fundamental group \cite{Morgan} \cite{HainKth}.  For simplicity we work 
with
the group algebra $\cc \pi _1^{\wedge }$ (completed at the augmentation
ideal).  The mixed Hodge structure on $\pi _1$ is usually defined via the mixed
Hodge structure on
 $\cc \pi
_1^{\wedge }$  and the inclusion of the Lie algebra corresponding to $\pi _1$
into this group algebra.

We first consider the relationship between the completed group algebra and the
completions of the spaces of representations at the identity.

Suppose $A$ is an augmented $\cc$-algebra which is complete with respect to
the augmentation ideal $J_A$.  Let $R(A,n)$ denote the functor of artinian 
local
$\cc$-algebras $B$ defined by setting 
$$
R(A,n)(B):= Hom ^{\rm aug}(A, M_n(B))
$$
where $Hom ^{\rm aug}$ denotes the set of algebra homomorphisms
sending $J_A$ to the ideal $M_n(\germ _B)$ ($\germ _B$ denotes the maximal
ideal of $B$).  

If $A=\cc \Gamma ^{\wedge }$ is the completion of the group algebra of a
finitely presented group $\Gamma$ then $R(A,n)$ is pro-representable by the
completion at the identity representation $R(\Gamma , GL(n))^{\wedge}$ of the
space of representations of $\Gamma$ in $GL(n)$ (there are probably more
abstract conditions on $A$ which could be used to insure 
representability but we don't need those here).  

Conversely let $\Cc$ denote the category of algebras which are direct 
produts of
algebras of the form $M_n(\cc )$. Suppose $\Upsilon : \Cc \rightarrow
\mbox{ForSch}$ is a functor from $\Cc$ to the category of formal schemes,
compatible with products (so we can think of $\Upsilon$ as a collection of
formal schemes $\Upsilon _n$ together with morphisms of functoriality
corresponding to morphisms of products of algebras in $\Cc$).  Then we define
${\cal A} (\Upsilon )$ to be the algebra of natural transformations $\Upsilon
\rightarrow 1_{\Cc }$.  The elements of ${\cal A}(\Upsilon )$ are functions $a$
which for each $n$ associate a section $a_n: \Upsilon _n\rightarrow M_n(\cc )$
with the $a_n$ compatible with morphisms of products of objects of $\Cc $.

\begin{lemma}
\label{Algs}
Suppose $A=\cc \Gamma ^{\wedge}$ is the completion of the group algebra of a
finitely presented group. Then the $R(A,n)$ give a functor $R(A,\cdot ):\Cc
\rightarrow \mbox{ForSch}$ and we can recover $A$ by the construction of the
previous paragraph:  
$$
A={\cal A}(R(A,\cdot )).
$$
\end{lemma}
\eop

{\em Remark:}
The morphisms of functoriality defining $R(A,\cdot )$ can be obtained from the
morphisms of functoriality of $R(\Gamma , GL(n))$ for morphisms between
products of groups $GL(n)$. 

Now we investigate the Hodge filtrations. Suppose $X$ is a smooth projective
variety. We obtain a family of formal completions $R_{Hod}(X,x,GL(n))^{\wedge}
\rightarrow {\bf A}^1$, with an action of $\Gm$.  The technique of Goldman
and Millson used in (\cite{Moduli} \S 10) to give an isomorphism
$R_{Dol}(X,x,GL(n))^{\wedge}\cong  R_{DR}(X,x,GL(n))^{\wedge}$ actually gives a
trivialization  $$
R_{Hod}(X,x,GL(n))^{\wedge}
\cong R_{Dol}(X,x,GL(n))^{\wedge}\times ^{\wedge} {\bf A}^1,
$$
with the action of $\Gm$ coming from the action defined in \cite{HBLS} on
$R_{Dol}(X,x,GL(n))^{\wedge}$ and the standard action on ${\bf A}^1$.
We discuss this further in \S 9 below.

We obtain a functor $R_{Hod}(X,x,\cdot ): \Cc \rightarrow \mbox{ForSch}/{\bf
A}^1$ which we think of as a family of functors parametrized by ${\bf A}^1$
(with action of $\Gm$). Because of the trivialization we can apply the
previous lemma.  This family of functors gives rise to a completed algebra
${\cal A}$ over ${\bf A}^1$, by a relative version of the construction of
Lemma \ref{Algs} (which poses no problem since everything is a product). 
Conversely starting from ${\cal A}$ we get back the
$R_{Hod}(X,x,GL(n))^{\wedge}$. Finally, the fiber ${\cal A}_1$ over $1\in {\bf
A}^1$ is isomorphic to the completed group algebra $\cc \pi _1(X,x)^{\wedge}$
again by the above discussion.  This family of algebras together with
$\Gm$-action (which as we have seen is equivalent to the data of the $R_{Hod}$
functorially in $n$) corresponds to a filtration of $\cc \pi
_1(X,x)^{\wedge}$.  We claim that this filtration is the Hodge filtration of
Morgan-Hain \cite{Morgan} \cite{HainKth}.

To see this, note a consequence of the trivializations and $\Gm$-actions
in the above discussion, that the filtration on $\cc \pi _1(X,x)^{\wedge}$
corresponding to ${\cal A}$ is just the filtration associated to the grading
given by the $\Gm$-action.  This $\Gm$-action is that which was defined in
\cite{HBLS}.  Finally, in \S\S 5-6 of \cite{HBLS} it was verified that this
$\Gm$-action gives rise to the Hodge filtration of Morgan-Hain.

To sum up, starting with the completed group algebra 
$\cc \pi _1(X,x)^{\wedge}$ and the Morgan-Hain Hodge filtration $F$ we can
form the family of algebras ${\cal A}=\xi ( \cc \pi _1(X,x)^{\wedge},F)$ over
${\bf A}^1$ with $\Gm$-action; then the family of completed representation
spaces associated to this family of algebras is isomorphic (together with
$\Gm$-action) to the completion $R_{Hod}(X,x,GL(n))^{\wedge}$ along the 
identity
section.  In the other direction, the data of the completions
$R_{Hod}(X,x,GL(n))^{\wedge}$ functorially in $GL(n)$ serve to define (via
Lemma \ref{Algs}) a family of completed algebras ${\cal A}$ over ${\bf A}^1$,
again with $\Gm$-action and isomorphism 
${\cal A}_1\cong \cc \pi _1(X,x)^{\wedge}$, and this family yields the
Morgan-Hain Hodge filtration on $\cc \pi _1(X,x)^{\wedge}$ by reversing the
construction $\xi$.

Thus the completion of our Hodge filtration at the identity representation
corresponds to the Hodge filtration on the nilpotent completion of the
fundamental group.  The whole Hodge filtration $R_{Hod}(X,x,G)$ or
$M_{Hod}(X,G)$ should be thought of as an analytic continuation of the Hodge
filtration on the nilpotent completion.

It might be interesting to try to express the existence of this analytic
continuation in terms of estimates on the mixed Hodge structure on 
$\cc \pi _1(X,x)^{\wedge}$ in the spirit of Hadamard's technique
\cite{Hadamard}. 

\subnumero{Writing formulas}

We can combine what we know so far to sketch a method which should allow,
in principle, to write down the local equations for the correspondence
$M_{Dol}\cong M_{DR}$ (and hence for the quaternionic structure) near a
complex variation of Hodge structure in $M^{\rm sm}$.  

The method sketched above should work equally well along any
prefered section coming from a complex variation of Hodge structure $\rho$.
The Hodge filtration on the relative Malcev completion \cite{HainMalcev} should
give a Hodge filtration $F$ on the complete local ring $\widehat{\Oo} _{M_{DR},
\rho}$; then taking $\xi (\widehat{\Oo} _{M_{DR},
\rho}, F)$ we get a family of complete local rings over ${\bf A}^1$, and
taking a formal spectrum we get a formal scheme along ${\bf A}^1$.  The real
structure or at least the invariant indefinite hermitian form underlying the
variation $\rho$ should give an involution allowing us to glue the formal
scheme with itself to get a formal scheme along $\pp ^1$. We should get back in
this way the formal completion of $M_{Del}$ along the prefered section. 

The
normal bundle of $M_{Del}$ along a prefered section is a direct sum of $\Oo
_{\pp ^1}(1)$'s, so the space of sections near the given section (which has a
structure of formal scheme in this case) will map isomorphically to the product
of any two fibers.  Taking two antipodal fibers (for example the fibers over
$0$ and $\infty$) we obtain explicitly the involution on the space of sections
and the prefered sections are those which are invariant.  Finally looking at
the isomorphism from the space of sections to the product of the formal
completions of $M_{Dol}$ and $M_{DR}$, the space of invariant sections gives
the graph of the real analytic isomorphism $M_{Dol}\cong M_{DR}$.  

One can imagine following out this entire construction explicitly to obtain
the Taylor series for the isomorphism $M_{Dol}\cong M_{DR}$ near the point
$\rho$.  The only ingredients are the Hodge filtration and the real structure
(and of course an analysis of the space of sections of our formal scheme, but
this is an algebraic question).

One can see just from the existence of this method that the algebraically
closed field generated by the coefficients of the Hodge filtration on the
relative Malcev completion (and their complex conjugates) will contain the
coefficients of the power series for the isomorphism $M_{Dol}\cong M_{DR}$ and
hence for the power series of the quaternionic structure.

I have not checked the details of this construction any more than what is
written above.

\numero{Formal categories}

One of the main properties of the Hodge filtration on usual abelian cohomology
is Griffiths transversality.  This is a property of the variation of the Hodge
filtration with respect to the Gauss-Manin connection which arises from a
smooth family of varieties.  We would like to obtain a similar property for
our nonabelian Hodge filtration.  Let's first look at how to interpret the
usual Griffiths transversality in terms of the construction $\xi$.

Suppose
$S$ is a smooth variety and $V$ is a vector bundle with integrable connection
$\nabla$.  Suppose $F^{\cdot}$ is a decreasing filtration of $V$ by subbundles.
Then $F^{\cdot}$ satisfies the Griffiths transversality condition $\nabla F^p
\subset F^{p-1}\otimes \Omega ^1_S$ if and only if the action of
$T(S)$ on $V$ extends to a $\Gm$-invariant action of the sheaf 
$T(S\times {\bf
A}^1/{\bf A}^1) (-S\times \{ 0\} )$ (of
relative tangent vector fields vanishing to order one along $S\times \{ 0\}$) 
on $\xi (V, F^{\cdot})$. This can be seen by calculating directly with the
definition of $\xi (V, F)$ (cf Lemma \ref{VBonXHod} below).  

In the nonabelian context suppose $X\rightarrow S$ is a smooth projective
morphism. We have a family $M_{DR}(X/S,G)$ of moduli spaces over $S$ and we
would like our ``Griffiths transversality'' to say that the lifting  of vector
fields on $S$ to vector fields on $M_{DR}(X/S,G)$ given by the Gauss-Manin
connection, extends to a $\Gm$-invariant lifting of sections of $T(S\times {\bf
A}^1/{\bf A}^1) (-S\times \{ 0\} )$ to vector fields on $M_{Hod}(X/S, G)$.

The difficulty in making this precise is that the Gauss-Manin
connection can no longer be interpreted in terms of vector fields if
$M_{DR}(X/S,G)$ is not smooth---so the above interpretation makes sense only on
the smooth points. The calculations in terms of vector fields are also 
difficult
to follow through.  To remedy these problems we introduce the point of view of
formal categories.  

Recall that the Gauss-Manin connection on $M_{DR}(X/S, G)$ is an isomorphism 
$$
p_1^{\ast}M_{DR}(X/S, G)| _{(S\times S)^{\wedge}} \cong 
p_2^{\ast}M_{DR}(X/S, G)| _{(S\times S)^{\wedge}}
$$
where $p_1, p_2: S\times S\rightarrow S$ are the projections and 
$(S\times S)^{\wedge}$ is the formal completion of the diagonal.  If we set
$N:=(S\times S)^{\wedge}$ then the pair $(S,N)$ has a structure of category in
the category of formal schemes.  The notion of formal category is a
generalization of this example. It provides a general framework for operations
on families of things over $S$.

A {\em formal category} is a pair $(X,N)$ consisting of a scheme $X$ and a
formal scheme $M$ mapping to $X\times X$, together with a structure of
category, that is morphisms $N\times _XN\rightarrow N$ giving composition and
$X\rightarrow N$ giving the identity, subject to the usual axioms for a
category. A formal category gives in a natural way a presheaf of
categories on $Sch/\cc $. A {\em formal groupoid} is a formal category such 
that
the values of the associated presheaf are groupoids.  We say that a formal
category is {\em of smooth type} if  $X$ is smooth,
the underlying scheme of $N$ is the scheme $X$ (via the identity morphism), and
$N$ is formally smooth. 

Let $X_N$ denote the stack over $Sch/\cc$ associated to the presheaf of
groupoids given by $(X,N)$.  We have a morphism $p:X\rightarrow X_N$. 
Note that $N$ represents the functor $X\times _{X_N}X$, so we can recover
$(X,N)$ from the stack $X_N$ with its morphism $X\rightarrow X_N$.  In
practice we confuse the notions (and notations) of formal groupoid $(X,N)$ and
associated stack $X_N$.

Suppose $(X,N)$ is a formal groupoid of smooth type. Note that the
structure sheaf $\Oo$ of $Scxh/\cc$ restricts to a sheaf of
rings which we also denote by $\Oo$ on $Sch /X_N$.  There is a complex of
$\Oo$-modules over $X_N$ which we denote by $p_{\ast}\Omega ^{\cdot}_{X/X_N}$
with differential denoted $d$, giving a resolution $$ \Oo \rightarrow 
p_{\ast}\Omega ^{0}_{X/X_N}\rightarrow  p_{\ast}\Omega ^{1}_{X/X_N}\rightarrow
\ldots \rightarrow p_{\ast}\Omega ^{n}_{X/X_N}\rightarrow 0
$$
(here $n=dim(N/X)$ is the relative (formal) dimension of the formally smooth
scheme $N$ over $X$ via either of the projections).

The notation $p_{\ast}\Omega ^{\cdot}_{X/X_N}$ is justified by
the fact that each component of this resolution is actually the direct image
of a locally free sheaf $\Omega ^i_{X/X_N}$ on $X$. This locally free sheaf
comes from $\Omega ^i_{N/X}$ by descent from $N=X\times _{X_N}X$ to $X$. Note
that the differential is in general a differential operator (of first order)
between these component sheaves, so it becomes a morphism only over $X_N$. 

A {\em local system $V$ on $X_N$}  is a sheaf on $Sch /X_N$ which locally
isomorphic to  $\Oo ^a$. This is equivalent to a vector bundle $V_X$ on $X$
together with {\em $N$-connection} $\varphi : V_X\rightarrow V_X \otimes _{\Oo
_X} \Omega ^1_{X/X_N}$ satisfying an integrability condition $\varphi ^2=0$.

We obtain
the resolution $p_{\ast}\Omega ^{\cdot}_{X/X_N}\otimes _{\Oo} V$ of $V$,
which we can use to calculate the cohomology of 
$V$ over $X_N$.  Note that the morphism $p$ is cohomologically trivial for
coherent sheaves on $X$. So the cohomology of $V$ may  be calculated as the
hypercohomology of the complex $\Omega ^{i}_{X/X_N}\otimes _\Oo V$ of Zariski
(or etale) sheaves on $X$; in particular there is a spectral sequence
$$
H^i(X, \Omega ^j_{X/X_N}\otimes _{\Oo}V)\Rightarrow H^{i+j}(X_N, V).
$$
The cohomology groups are finite dimensional vector spaces if $X$ is
proper.

To a formal category $(X,N)$ of smooth type we can associate a split almost
polynomial sheaf of rings of differential operators $\Lambda _N$
\cite{Illusie} \cite{Moduli}.  It is the sheaf of rings associated to the
differentials in the above complex. Note that $p_{\ast}(\Oo _N )$ is naturally
a projective limit of locally free sheaves on $X$.  We can construct $\Lambda
_N$ as the continuous dual, which is a union of locally free sheaves. The ring
structure is dual to the cogebra structure
$$
p_{\ast}(\Oo _N )\rightarrow p_{\ast}(\Oo _N )\otimes _{\Oo _X}p_{\ast}(\Oo _N
) 
$$
which itself comes from the composition morphism $N \times _XN\rightarrow N$.
A local
system $V$ on $X_N$ is the same thing as a $\Lambda _N$-module; the underlying
$\Oo_X$-module is $V_X$ and the $\Lambda _N$-module structure is given by
$\varphi$.  

A {\em principal $G$-bundle} on $X_N$ means a $G$-torsor on the stack $X_N$. 
This is the same thing as a principal $G$-bundle $P$ on $X$ together with
an isomorphism $\varphi :s^{\ast}(P)\cong b^{\ast}(P)$ on $N$ (where
$s,b:N\rightarrow X$ are the two tautological morphisms), such that
$\varphi$ satisfies the appropriate cocycle condition on $N\times _XN$.

There is a notion of semistability for local systems over $X_N$ which is
analogous to the usual notion: we can define a notion of coherent sheaf over
$X_N$ (a coherent sheaf on $X$ with descent data to $X_N$) and a bundle over
$X_N$ is semistable if for every $X_N$-subsheaf, the normalized Hilbert
polynomial is less than or equal to that of the original object. Again as usual
we can define the notion of semistability of a principal $G$-bundle on $X_N$.
Lacking a direct proof of the conservation of semistability by tensor product
in the case of local systems over a general formal category (this is a good
question for further research), we put in the definition here that the local
systems associated to all representations of $G$ should be semistable.

When considering a relative situation $X_N\rightarrow S$, semistability means
semistability on each fiber $(X_N)_s$ (it is an open condition on the 
base $S$).

Because, in all of the examples we interested in in this paper, it is
necessary to include a condition of vanishing rational Chern classes when
defining the moduli spaces, we put this directly into the definition in the
formal category setting.  Of course, for formal groupoids different from our
examples, this condition may not necessarily be a sensible one; and even in our
examples, it may also be interesting to consider other components of the
moduli spaces.  One could make the same definitions and obtain moduli spaces
without this condition, but we include the condition here for simplicity of
notation.

Suppose $(X,N)\rightarrow S$ is a morphism from a formal groupoid of smooth
type to a base scheme $S$, such that $X$ is smooth and proper over $S$ and $N$
is formally smooth over $S$. Suppose $x: S\rightarrow X$ is a section. Define
the functor $\Rr (X_N/S, x, G)$ which to an $S$-scheme $S'$ associates $\{
(P,\varphi , \beta )\}$ where $(P, \varphi )$ is a principal
$G$-bundle over $X_N\times _SS'$, semistable and with vanishing rational Chern
classes relative to $S$';  and $\beta : x^{\ast}(P) \cong G\times S'$ is a
framing along the section $x$.

\begin{theorem}
\label{ConstructMN}
The functor $\Rr (X_N/S, x,
G)$ is representable by a scheme which we denote by $R(X_N/S, x, G)\rightarrow
S$. Furthermore all points are semistable for the action of $G$ so a universal
categorical quotient $M(X_N/S,G)= R(X_N/S,x, G)// G$ exists.   
\end{theorem}
This theorem follows from the interpretation of local systems over $X_N$ as
$\Lambda _N$-modules \cite{Moduli}; the construction of the moduli and
representation schemes for $\Lambda _N$-modules; and from the tannakian point
of view used in \cite{Moduli} \S 9.
\eop

We denote the {\em stack quotient} by 
$$
\Mm (X_N/S,G):=  R(X_N/S,x, G)/ G.
$$ This is the first
relative nonabelian cohomology of $X_N/S$ with coefficients in $G$.
It is an algebraic stack.

\subnumero{The basic examples}

The main example (which we used to introduce this section and which you meet in
most treatments of the subject---cf \cite{Berthelot} for example) is the formal
groupoid obtained by setting $N:= (X\times X)^{\wedge}$ (completion along the
diagonal).  We denote this formal groupoid (or the associated stack) by
$X_{DR}$.  In this stack there is at most one morphism between any pair of
objects, so the stack is equivalent to a sheaf of sets which we also denote
$X_{DR}$.  The sheaf of sets is the quotient of $X$ by the equivalence relation
$N$: heuristically we identify any two points which are infinitesimallly close
together, and it is this infinitesimal glue which makes it so that $X_{DR}$
actually reflects the topology of the underlying usual space.  More precisely
if $S$ is any scheme over $\cc$ then the $S$-valued points of $X_{DR}$ are
the $S$-valued points of $X$ modulo the relation that two
points are equivalent if their restrictions to the underlying reduced scheme
$S^{\rm red}$ are the same; except that we have to divide by this
equivalence relation and then sheafify.  Since $X$ is smooth, any $S^{\rm
red}$-valued point extends, locally on $S$, to an $S$-valued point.  Thus after
taking the quotient and sheafifying the result is simply that
$X_{DR}(S)=X(S^{\rm red})$.  

The sheaf of rings of differential operators associated to the formal groupoid
$X_{DR}$ is just the full ring $\Lambda _{DR}$ of differential operators on
$X$ \cite{Moduli}.  A principal bundle on $X_{DR}$ is
just a principal bundle on $X$ with integrable connection, and a vector bundle
or local system over $X_{DR}$ is just a vector bundle on $X$ with integrable
connection. 
We recover 
$$
M(X_{DR}, G)=M_{DR}(X,G)
$$
and similarly for the representation spaces and stacks 
$$
R(X_{DR},x, G)=R_{DR}(X,x,G), \;\;\;\; \Mm (X_{DR}, G)=\Mm _{DR}(X,G).
$$
The cohomology of $X_{DR}$ with coefficients in a local system is just the
algebraic de Rham cohomology of $X$ with coefficients in the corresponding
vector bundle with integrable connection.

We now define a formal groupoid $X_{Dol}$ which gives rise to the Dolbeault
theory in the same way as $X_{DR}$ gave rise to the de Rham theory.  
In this formal groupoid 
the object object is $X$ and the morphism object is the formal completion of 
the
zero section in the tangent bundle of $X$, lying over the diagonal in $X\times
X$. A principal bundle on $X_{Dol}$ is just a principal Higgs bundle; a local
system is  a Higgs bundle; the cohomology of a local system is the Dolbeault
cohomology; and the associated sheaf of rings of differential operators is the
ring $\Lambda _{Dol}$ defined in \cite{Moduli}.  
We recover 
$$
M(X_{Dol}, G)=M_{Dol}(X,G)
$$
and similarly for the representation spaces and stacks.

Now we define a formal groupoid $X_{Hod} \rightarrow {\bf A}^1$ which serves
as a deformation from $X_{DR}$ to $X_{Dol}$ and from which we can recover the
moduli spaces $M_{Hod} (X, G)$. It is the stack associated to the
realisation of the nerve of the presheaf of groupoids given by a formal
groupoid which we denote by $\widetilde{X}_{Hod} \rightarrow {\bf A}^1$. The
object object is   $$ Ob(\widetilde{X}_{Hod}):=X\times {\bf A}^1.
$$
Let $Y$ be the complement of the strict transform of $X\times X\times \{ 0\}$ 
in
the blow-up of $X\times  X\times {\bf A}^1$ along $\Delta (X) \times \{ 0\}$.
(Here $\Delta (X)$ is the diagonal.)  There is a unique composition
$$
Y\times _{X\times {\bf A}^1}Y \rightarrow Y
$$
compatible with the trivial composition 
$$
(X\times X\times {\bf A}^1)\times _{X\times  {\bf A}^1}
(X\times X\times {\bf A}^1)
$$
via the morphism $Y\rightarrow X\times X\times {\bf A}^1$.  

There is a morphism $\Delta ':X \times {\bf A}^1 \rightarrow Y$ covering the
inclusion $\Delta : X\rightarrow X\times X\times {\bf A}^1$. This section
provides an identity for the above composition---in particular, setting the
morphism object equal to $Y$ would define a category. As with the case of
$X_{DR}$ itself, we take the formal completion of this morphism object: the
morphism object $Mor(  \widetilde{X}_{Hod })$ is defined to be the formal
completion of $Y$ along  $\Delta '(X\times {\bf A}^1)$.

The formal groupoid defined in this way
is a groupoid;
it maps to ${\bf A}^1$;
it has fiber over $\{ 0\}$ equal to the formal groupoid defining $X_{Dol}$;
and its fiber over $\{ t\}$ for any $t\neq 0$ is equal to the formal groupoid
defining $X_{DR}$.

The general construction of Theorem \ref{ConstructMN} gives back for $X_{Hod}$
over ${\bf A}^1$ the result of Proposition \ref{ConstructMhod}:
$$
M(X_{Hod}/{\bf A}^1, G)= M_{Hod}(X,G)
$$
and similarly for the representation space and moduli stack
$$
R(X_{Hod}/{\bf A}^1,x, G)=R_{Hod}(X,x,G), \;\;\;\; \Mm (X_{Hod}/{\bf A}^1,
G)=\Mm _{Hod}(X,G). 
$$ 

\begin{lemma}
\label{VBonXHod}
In the case $G=GL(n)$ a section of the moduli stack 
${\bf A}^1\rightarrow \Mm (X_{Hod}/{\bf A}^1,GL(n))$ preserved by $\Gm$ (or
more precisely with action of $\Gm$ specified) corresponds to a
vector bundle with filtration  satisfying Griffiths transversality. The
relative cohomology of such a family of local systems is just the sheaf  over
${\bf A}^1$ with action of $\Gm$ corresponding to the induced filtration on the
cohomology of the local system.  
\end{lemma}
{\em Proof:}
A section of the moduli stack with action of $\Gm$ is just a vector bundle on
$X_{Hod}$ with action of $\Gm$.  In particular we have a vector bundle on
the underlying scheme
$X\times {\bf A}^1$ together with action of $\Gm$; by a relative version of
the inverse of $\xi$ this corresponds to a bundle $V$
on $X$ with filtration by subbundles $F^p$ such that the associated-graded is a
bundle. The descent data to $X_{Hod}$ are determined by the descent data over
$\Gm \subset {\bf A}^1$, and by the $\Gm$-invariance of our section, these
are determined by the descent data over $1\in {\bf A}^1$, which is to say an
integrable connection on $V$.  The original bundle over $X\times {\bf A}^1$  
is the locally free sheaf $\xi (V,F)=\sum t^{-p} F^p$.  The statement that the
connection extends to descent data for this bundle down to $X_{Hod}$ is
equivalent to the condition 
$$
t\nabla (\sum t^{-p} F^p)\subset (\sum t^{-p}
F^p)\otimes \Omega ^1_X,
$$
which translates to $\nabla F^p \subset
F^{p-1}\otimes \Omega ^1_X$---the Griffiths transversality condition.
\eop

\numero{The Gauss-Manin connection}
Suppose $X\rightarrow S$ is a smooth projective morphism.  Then we obtain
$$
M_{DR}(X/S, G)\rightarrow S,
$$
a family whose fiber over $s\in S$ is $M_{DR}(X_s,G)$. This family has an
algebraic integrable connection \cite{NAHT} \cite{Moduli}.  We can interpret 
the
connection in terms of formal categories in the following way (this is a
simple variant of the crystalline interpretation of \cite{Moduli}). We have a
morphism $X_{DR}\rightarrow S_{DR}$ and the fiber product
$X_{DR}\times _{S_{DR}}S$ has a structure of smooth formal groupoid over $S$,
which we call $X_{DR/S}$.    The morphism space is the formal completion of
the diagonal in $X\times _SX$.  The moduli stack $\Mm _{DR}(X/S)$ is just the 
nonabelian cohomology of $X_{DR/S}$ relative to $S$ with coefficients in $G$,
or in our previous notations $$
\Mm _{DR}(X/S, G) = \Mm (X_{DR/S}/S, G).
$$  
The same holds for the moduli spaces and representation spaces.
But we could equally well take the
nonabelian cohomology of $X_{DR}$ relative to $S_{DR}$.  We obtain a stack over
$S_{DR}$ which, when pulled back to $S$, gives $\Mm _{DR}(X/S, G)$.  To put 
this
another way, we get descent data for $\Mm _{DR}(X/S, G)$ from $S$ down to
$S_{DR}$.  This is exactly the data of an integrable connection which is
the nonabelian version of the Gauss-Manin connection.    

One can see that the associated analytic connection on the analytic family is
the same as that induced by the fact that (locally over the base) all of the
fibers of $\Mm _{DR}(X/S)$ are of the form $\Mm _B(\Gamma )$ for $\Gamma $ the
fundamental group of the fiber (\cite{Moduli}, Theorem 8.6).

In abelian Hodge theory there are two principal results about the Gauss-Manin
connection: Griffiths transversality with respect to the Hodge filtration, and
regular singularities at the singular points of a family. We obtain their
analogues for nonabelian cohomology by using the theory of formal categories
and following the above description of the connection.

These properties are easily obtained by using a variant of our construction of
the connection.  Suppose $X_N \rightarrow S_K$ is a morphism of formal
categories of smooth type such that the fiber product $X_N\times _{S_K} S$ is a
formal groupoid of smooth type on $X/S$.  Then the schemes $M(X_N\times
_{S_K}S/S,G)$ and stacks  $\Mm (X_N\times _{S_K}S/S,G)$ have descent data down
to $S_K$, that is they are pullbacks of sheaves or stacks on $Sch /S_K$.  The
same is true for the  $R(X_N\times _{S_K}S/S,x, G)$ if $x: S_K\rightarrow X_N$
is a section. 

\subnumero{Griffiths transversality}

Suppose $X\rightarrow S$ is a smooth family of projective varieties (and we
now ask that the base be smooth, although we don't need it to be
projective).  Then we obtain a morphism of formal categories 
$X_{Hod} \rightarrow S_{Hod}$ over ${\bf A}^1$.  Put
$$
X_{Hod /S} := X_{Hod} \times _{S_{Hod}} (S\times {\bf A}^1).
$$
It is given by a smooth formal groupoid on $X\times {\bf A}^1$
relative to $S\times {\bf A}^1$. The relative
nonabelian cohomology is
$$
\Mm _{Hod}(X/S, G):= \Mm (X_{Hod/S}/S, G)\rightarrow S\times {\bf A}^1.
$$  
This  morphism is provided with an action of the formal groupoid $S_{Hod}$
(i.e. $\Mm _{Hod}(X/S, G)$ is the pullback of  a stack over $S_{Hod}$).  In
particular over $\lambda \neq 0$ we get back the action of $S_{DR}$, that is to
say the Gauss-Manin connection.  There is a $\Gm$-action
compatible with everything. This whole situation is the nonabelian
analogue of Griffiths transversality, an interpretation which, comparing with
the abelian case, is justified by Lemma {VBonXHod}.

After we discuss the compactification below we will give another
interpretation in terms of poles of the Gauss-Manin connection at infinity.

\subnumero{Regularity of the Gauss-Manin connection}

Suppose $S'\subset S$ is an open subset whose complement is a
divisor $D$ with normal crossings.  Define $S_{DR}(\log D)$ to be the
formal groupoid which is associated to
the bundle of vector fields tangent to $D$.  To construct it
explicitly we first treat the case where $S$ is a smooth curve and $D$ a point. 
Then $S\times S$ is a smooth surface.  Blow up at the point $(D,D)$,  
and take the formal completion along
strict transform of the diagonal, as
morphism scheme.  The maps of a smooth scheme into the blow up minus the
transform of $D\times S$ are the same as maps into $S\times
S$ which send the intersection with the divisor $D\times S$ to the
point $(D,D)$.  From this description we obtain the composition. The formal
completion then becomes a formal groupoid.  Now, for any $(S, D)$, glue this
construction in along each component of the divisor (or more precisely along
an $i$-fold intersection of the divisor in an $n$-dimensional $S$,
glue the product of $i$ copies of this construction with $n-i$ copies of
a smooth curve with the usual de Rham construction.

A vector bundle or local system over $S_{DR} (\log D)$ is just a vector bundle
on $S$ with integrable connection with logarithmic singularities along $D$. In
particular, the condition that a vector bundle on $S'_{DR}$
(corresponding to a vector bundle with integrable connection on $S'$) extends 
to
a bundle on $S_{DR}(\log D)$ is equivalent to the condition that the
connection have regular singularities. 

Now suppose that we have a projective morphism of smooth
varieties $f:S\rightarrow S$ and a divisor $D\subset
S$ which has  normal crossings, such that  $Y:=f^{-1}(D)$ has normal
crossings, and such that $f:X'\rightarrow S'$ is smooth where $S'=S-D$
and $X'=f^{-1}(S')=X-Y$.  
There is a morphism of formal groupoids $S_{DR}(\log D)\rightarrow S_{DR}$.  Let
$$
X_{DR} (\log D):= X_{DR}\times _{S_{DR}}S_{DR}(\log D).
$$
Again we have a morphism of formal groupoids
$$
X_{DR} (\log D)\rightarrow S_{DR} (\log D),
$$
and the fiber product 
$$
X_{DR/S} (\log D):= X_{DR} (\log D)\times _{S_{DR} (\log D)}S
$$
is a formal groupoid on $X$ over $S$. Even though $X/S$ is not smooth, this
formal groupoid corresponds to a split almost polynomial sheaf of rings of
differential operators $\Lambda _{DR/S (\log D)}$ on $X$ over $S$.  Thus the
moduli problems are solved by \cite{Moduli} so the first relative nonabelian
cohomology stack $\Mm (X_{DR/S} (\log D)/S, G)$ is an algebraic stack; the
associated representation functor for objects provided with a frame along a
fixed section (not passing through the singular points) is
representable by a scheme, and the universal categorical quotient scheme
$M(X_{DR/S} (\log D)/S,G)$ exists. 

Finally, since 
$X_{DR/S}(\log D)$ comes from fiber product as in our general situation, these
moduli spaces, moduli stacks and representation spaces descend to $S_{DR}(\log
D)$. This statement is the regularity of the Gauss-Manin connection.  

If $G=GL(n)$ then we can do slightly better. Recall that   
$M(X_{DR/S} (\log D)/S,GL(n))$ 
parametrizes bundles with descent data to the formal groupoid $X_{DR/S} (\log
D)$, which are semistable with vanishing rational Chern classes on the fibers
$X_s$.   But over the singular fibers, where the extra structure is a
connection away from the singularities but simply a logarithmic connection near
the singularities, it is also natural to consider objects which are no longer
bundles but torsion-free sheaves.  This moduli space, which we can denote by
$M^{\rm tf}(X_{DR/S} (\log D)/S,GL(n))$, 
will have the advantage that, when combined with the techniques used below for
compactifying $M_{DR}(X_s, GL(n))$, will give a compact total space over
projective $S$.

{\em Caution:} One might be tempted to try the above argument with
$X_{DR}/S_{DR}$.  In this case the fiber product $X_{DR}\times _{S_{DR}}S$
no longer has the required smoothness properties, and the relative moduli stack
is no longer an algebraic stack---otherwise we would have an extension of the
Gauss-Manin connection over the singularities! In the case $G= {\bf G}_a$ for
example this would give an extension of the Gauss-Manin connection for ordinary
first cohomology, easily seen as impossible in most examples.

{\em Exercise:}  Explain what goes wrong in an example (such as a family of
curves aquiring a node) if we try to do the previous construction with
$X_{DR}/S_{DR}$.  This can be done in the context of abelian cohomology.

\numero{Etale local triviality of $M_{Hod}$}

\subnumero{Goldman-Millson theory}
A deformation problem is often
controlled by a differential graded Lie algebra (dgla) $L^{\cdot}$ over the
base field which we are assuming is $\cc$.  As explained in
\cite{GoldmanMillson} this means that the deformations with values in an artin
local ring $A$ (with maximal ideal $\germ$) correspond to elements $\eta \in
L^1\otimes _{\cc} \germ$ such that  $$
d(\eta ) + \frac{1}{2}
[\eta , \eta ] = 0.
$$
The isomorphisms between deformations $\eta $ and $\eta '$ correspond to
elements $s\in L^0\otimes _{\cc} \germ$ with 
$$
\eta ' = d(s) + e^{-s} \eta e^s .
$$

If $P$ is a principal $G$-bundle, let $A^{\cdot} (ad \, P)$ denote the graded
Lie algebra of $\Cc ^{\infty}$ forms on $X$ with coefficients in the adjoint
bundle $P\times ^G\gerg$ (with Lie bracket combining the wedge of forms and the
Lie bracket of $\gerg$). If $(P,\theta )$ is a principal Higgs bundle structure
then we obtain a  dgla $(A^{\cdot}(ad\, P), \delbar + \theta )$ which
gives the deformation theory of $(P,\theta )$. If $(P,\nabla )$ is a principal
bundle with integrable connection then we obtain 
a  dgla $(A^{\cdot}(ad\, P), \nabla  )$  which gives
the deformation theory of $(P,\nabla  )$ \cite{GoldmanMillson}.

More generally, suppose $(X,N )$ is a formal groupoid structure for
$X$; then we have the algebra of differentials 
$\Omega ^{\cdot}_{X/X_N}$ (cf \S 7 above). Let 
$$
A^{i}_{N}:= \bigoplus _{p+q=i} A^{0,q} (\Omega ^p_{X/X_{N}}
$$
This is a differential graded algebra with differential equal to $\delbar +
\delta $ where $\delta$ is the first order differential operator corresponding
to the differential of $p_{\ast}\Omega ^{\cdot}_{X/X_{N}}$. 
If $P$ is a principal $G$-bundle over $X_{N}$ then we obtain a dgla
$$
A^{\cdot}_{N}(ad\, P):= A^{\cdot}_{N} \otimes _{\Oo} (P\times ^G\gerg ).
$$
The differential comes from that of $A^{\cdot}_{N}$. This dgla controls the
deformation theory of $P$ as a principal bundle on $X_{N}$ (we leave the proof
as an exercise following \cite{GoldmanMillson} and \cite{Moduli} \S 10). 

If
$X_{N} = X_{DR}$ then we get $A^{\cdot}_{DR}(ad\, P):=(A^{\cdot}(ad\, P),
\nabla  )$, whereas if $X_{N} = X_{Dol}$ then we get  $A^{\cdot}_{Dol}(ad\,
P):=(A^{\cdot}(ad\, P), \delbar + \theta )$.

If $L^{\cdot}$ is a dgla then let $H^{\cdot}(L^{\cdot})$ denote the dgla of
cohomology with differential equal to zero.  Recall that we say that
$L^{\cdot}$ is {\em formal} if there is a quasiisomorphism between
$L^{\cdot}$ and $H^{\cdot}(L^{\cdot})$. According to the theory of Goldman and
Millson \cite{GoldmanMillson}, a quasiisomorphism induces an equivalence of
deformation theories, and on the other hand the deformation theory of a dgla
with zero differential is {\em quadratic}, in other words the universal
deformation space is the quadratic cone in $H^1$ which is defined as the zero
scheme of the map $H^1 \rightarrow H^2$, $\eta \mapsto [\eta , \eta ]$.  It
follows that the deformation theory of a formal dgla is quadratic.  

We need a relative version of this theory for $X\times {\bf A}^1/{\bf A}^1$
near a prefered section.   It is actually an interesting question (which I 
don't
think has yet been addressed) to develop a relative version of
Goldman-Millson theory in all generality.  In our case we are helped by the
fact that the total space is a product.  Also we will restrict to deformations
over artinian base (whereas in a better version one should consider
arbitrary base with nilpotent ideal).  

Suppose $P$ is a flat principal $G$-bundle with harmonic $K$-reduction $P_K$
and associated operators $d'=\partial + \theta '$ and $d'' = \delbar + \theta
''$.  We define the differential graded Lie algebra over $\cc [\lambda ]$
$$
A^{\cdot}_{Hod}:=(A^{\cdot} (ad \, P) \otimes _{\cc} \cc [\lambda ],
\lambda \partial + \theta ' + \delbar + \lambda \theta '' ).
$$
Note that the differential has square zero and, as varying with parameters,
gives a deformation from the de Rham to the Dolbeault differentials.  
Notice also that the components of $A^{\cdot}_{Hod}$ are flat $\cc [\lambda
]$-modules.  For any dgla $L^{\cdot}$ over $\cc [\lambda ]$ where the
components are flat, we define a stack of groupoids over the category of
artinian local $\cc [\lambda ]$-algebras $(B, \germ )$ in the same way as
above (except that tensor products are taken over $\cc [\lambda ]$): the
objects are elements  
$\eta \in L^1\otimes _{\cc [\lambda ]} \germ$ such that 
$$
d(\eta ) + \frac{1}{2}
[\eta , \eta ] = 0.
$$
The isomorphisms between deformations $\eta $ and $\eta '$ correspond to
elements $s\in L^0\otimes _{\cc [\lambda ]} \germ$ with 
$$
\eta ' = d(s) + e^{-s} \eta e^s .
$$
Again we have the same theorem that quasiisomorphisms between flat dgla's
induce equivalences of deformation groupoids. And, on the other hand, the dgla
$A^{\cdot}_{Hod}$ gives the deformation theory of $\Mm _{Hod}$ near the
prefered section corresponding to $P$.  More precisely the groupoid defined
above for $(B,\germ )$ is equivalent to the groupoid of morphisms
$Spec (B)\rightarrow \Mm _{Hod}$ with isomorphism between the morphism $Spec
(B/\germ )\rightarrow \Mm _{Hod}$ and the composed morphism 
$Spec
(B/\germ )\rightarrow {\bf A}^1 \rightarrow \Mm _{Hod}$ (where the first
morphism is the projection to ${\bf A}^1=Spec (\cc [\lambda ])$ and the second
morphism is the prefered section corresponding to $P$).  The proofs are
left as exercises.

Finally, the dgla $A^{\cdot}_{Hod}$ defined above is formal over $\cc [\lambda
]$.  Let $\ker(\partial + \theta '')[\lambda ]$ denote the subcomplex of
$A^{\cdot}_{Hod}$ consisting of forms $\alpha$ with $\partial \alpha + \theta
'' \alpha = 0$. Since this operator doesn't depend on $\lambda$, the subcomplex
is just a tensor product of the usual complex $\ker(\partial + \theta '')$ with
$\cc [\lambda ]$. Furthermore the subcomplex is a sub-dgla; and finally note
that the differential of the sub-complex is $\delbar + \theta '$ which also
doesn't depend on $\lambda$. The ``principle of two types''  
(cf Lemma 2.2 of \cite{HBLS}) implies that
the morphisms  
$$ 
H^{\cdot} (A^{\cdot}_{Hod}) \leftarrow \ker (\partial + \theta
'' )[\lambda ]  \rightarrow A^{\cdot}_{Hod}
$$
are quasiisomorphisms. 
Finally, let $\Hh ^{\cdot}$ denote the $\cc$-vector space of harmonic forms in
$A^{\cdot}(ad \, P)$ (which is the same for the flat connection or the
Dolbeault operator $\delbar + \theta '$ or, for that matter, anything in
between). Then the morphism  $$ \Hh ^{\cdot } \otimes _{\cc } \cc [\lambda ]
\rightarrow  H^{\cdot} (A^{\cdot}_{Hod})  $$ is an isomorphism of graded vector
spaces (the space of harmonic forms does not {\em a priori} have a product
structure).  This shows that  $H^{\cdot} (A^{\cdot}_{Hod})$ is flat over $\cc
[\lambda ]$ so we can apply the quasiisomorphism result to conclude that the
deformation theory of $A^{\cdot}_{Hod}$ is the same as that of the graded Lie
algebra  $H^{\cdot} (A^{\cdot}_{Hod})$.  Finally, in order to obtain
local triviality we have to show that there is a product structure on $\Hh
^{\cdot }$ such that the isomorphism 
$H^{\cdot} (A^{\cdot}_{Hod})\cong \Hh ^{\cdot } \otimes _{\cc } \cc [\lambda ]
$ becomes an isomorphism with product structure.  This can be seen by
identifying $\Hh ^{\cdot}$ with the cohomology of the complex 
$\ker(\partial + \theta '')$ and noting that this latter has a product
structure.  With this result, the deformation theory
along the prefered section becomes a product.  We have shown this result for
deformations over artinian base, so we get the result on the level of formal
completions.  Artin approximation then gives that locally in the etale topology
at any point of a prefered section, $\Mm _{Hod}$ is a product. Any point (even
non-semisimple) is isomorphic to a point in a neighborhood of a semisimple
point, so we obtain local triviality at any point in the union $\Mm
_{Hod}(X,G)$ of components corresponding to bundles with vanishing rational
Chern classes.

The formal trivialization along the prefered section is the total-space version
of the {\em isosingularity principle} stated in the introduction and
\S 10 of \cite{Moduli}.

We can sum up in the following theorem.

\begin{theorem}
\label{LocTriv}
Suppose $X$ is a smooth projective variety.  Let
$M_{Hod}(X,G)\rightarrow {\bf A}^1$ denote the union of components
corresponding to objects with vanishing rational Chern classes. Then etale
locally (above) $M_{Hod}(X,G)$ is a product, in other words any point $P\in
M_{Hod}(X,G)$ over $\lambda \in {\bf A}^1$ admits an etale neighborhood $P\in
U\rightarrow M_{Hod}(X,G)$ with an etale morphism $U\rightarrow
M_{Hod}(X,G)_{\lambda }\times  {\bf A}^1$. The same holds for $R_{Hod}
(X,x,G)$ and the moduli stack $\Mm _{Hod}(X,G)$. \end{theorem} \eop

\begin{corollary}
The morphism $M_{Hod}(X, G, 0) \rightarrow {\bf A}^1$ is flat (and similarly
for  $R_{Hod} (X,x,G)$ and the moduli stack $\Mm _{Hod}(X,G)$).
\end{corollary}
\eop

On the other hand, if we have a family of varieties $X\rightarrow S$ then 
the Gauss-Manin connection guarantees that the family $M_{DR}(X/S,G)\rightarrow
S$ is etale locally a product.  It seems reasonable to guess that these two
results combine into the following.

\begin{conjecture}
\label{triv}
Suppose $X\rightarrow S$ is a smooth projective family.  Let
$M_{Hod}(X/S,G)\rightarrow S\times {\bf A}^1$ denote the relative $M_{Hod}$
space. Then etale locally (above) $M_{Hod}(X/S,G)$ is a product, in other 
words any
point $P\in M_{Hod}(X/S,G)$ over $(s,\lambda )\in S\times {\bf A}^1$ admits an
etale neighborhood $P\in U\rightarrow M_{Hod}(X/S,G)$ with an etale morphism
$U\rightarrow M_{Hod}(X/S,G)_{(s,\lambda )}\times S\times {\bf A}^1$.
The same holds for $R_{Hod} (X/S,x,G)$ and the moduli stack 
$\Mm _{Hod}(X/S,G)$.
\end{conjecture}

To prove this one would have to analyze much more closely the deformation
theory. In particular, the fact that the total space is a product, which
helped a lot in the previous argument, is no longer there to help us here.

We do give an argument showing this conjecture over smooth points (i.e.
showing that the morphism to $S$ is smooth at smooth points of the fibers) in
the subsection ``Griffiths transversality revisited'' of \S 11 below.

A consequence of this conjecture would be that the deformation class of
$M_{Hod}$ over a point $(s, 0)$ (we take $\lambda = 0$ because the theory over
$\lambda \neq 0$ is trivial due to the Gauss-Manin connection) is determined by
a class  $\zeta \in H^1(M_{Dol}(X_s), \Theta_{M_{Dol}})\otimes (TS_s \oplus \cc
)$, where $\Theta_{M_{Dol}}$ is the etale sheaf of infinitesimal automorphisms
of  $M_{Dol}$.

\numero{A weight property for the Hodge filtration}

From its very definition, the usual Hodge filtration on cohomology has the
property that $F^0 V=V$. In terms of the Rees bundle this translates to the
statement that if $v\in V$ then the $\Gm$-orbit $\Gm v$ has a limit point,
i.e. it extends to a section ${\bf A}^1 \rightarrow \xi (V,F )$.  We will
establish a similar property for $M_{Hod} $ and ${\cal M}_{Hod}$. 

We don't explicitly review the notion of sheaf of rings of differential
operators $\Lambda$ on $X/S$, for $X\rightarrow S$ a projective morphism, from
\cite{Moduli}. Recall that for any formal groupoid of
smooth type we obtain an almost-polynomial sheaf of rings of differential
operators  as described in \cite{Illusie} and in \S 7 above.  The
$\Lambda$-modules are just the coherent sheaves on $X$ with descent data down
to $X_M$. 

On the other hand, in the case which interests us (the formal groupoid
$X_{Hod}$ on $X\times {\bf A}^1$ over ${\bf A}^1$) we can explicitly describe
the sheaf of rings $\Lambda$ (it is the sheaf of rings denoted $\Lambda ^R$
in \cite{Moduli}).  Note first of all that $\Lambda _{DR}$ (corresponding to
$X_{DR}$) is just the sheaf of rings of differential operators.  It has a
filtration $\Lambda _{DR}^i$ being the differential operators of order $\leq
i$.  Define a decreasing filtration by indexing negatively, $F^{-i}=\Lambda
^i_{DR}$.  This filtration is compatible with the ring structure so the
construction $\xi$ gives  a sheaf of rings on $X\times {\bf A}^1$ over ${\bf
A}^1$,
$$
\Lambda _{Hod} = \xi (\Lambda _{DR},F).
$$
It is the sheaf of rings associated to the formal groupoid $X_{Hod}$.
The relative moduli theory for $\Lambda _{Hod}$-modules on $X\times {\bf
A}^1/{\bf A}^1$ yields the moduli space $M_{Hod}$ and   representation space
$R_{Hod}$. The stack theoretic quotient gives the moduli stack $\Mm _{Hod}$.

\subnumero{Langton theory} 

We recall the notations and terminology of \cite{Moduli}.  In particular,
$p$-semistablity and $p$-stability refer to Gieseker's definition involving
Hilbert polynomials. We will work with $G=GL(n)$ at the start.

The following theorem is the generalisation of Langton's theory
\cite{Langton} of properness of moduli spaces, in the case of sheaves of
$\Lambda$-modules (M. Maruyama pointed out to me that Langton's theory carries
over in this type of general context).

\begin{theorem}
\label{ThmA}
Suppose $S=Spec (A)$ where $A$ is a discrete valuation ring with fraction field
$K$ and residue field $A/{\bf m} = \cc$. Let $\eta$ denote the generic
point and $s$ the closed point of $S$. Suppose $X\rightarrow S$ is a projective
morphism of schemes with relatively very ample $\Oo (1)$ on $X$.  Suppose
$\Lambda$ is a split almost polynomial sheaf of rings of differential operators
on $X/S$  as in \cite{Moduli}. Suppose $\Ff $ is a sheaf of
$\Lambda$-modules on $X$ which is relatively of pure dimension $d$, flat over
$S$, and such that the generic fiber $\Ff _{\eta}$ is $p$-semistable.  Then
there exists a sheaf of $\Lambda$-modules $\Ff '$ on $X$ which is relatively of
pure dimension $d$, flat over $S$, and such that $\Ff '_{\eta} \cong \Ff
_{\eta}$ and also $\Ff '_s$ is $p$-semistable. 
\end{theorem}
{\em Proof:}
Langton's proof carries over into our situation.  We give a brief sketch 
for compatibility with our notations. Let $p_X(\cdot )$ denote the absolute
normalized Hilbert polynomial of a sheaf on $X$ with proper support, and let
$p_{X/S}(\cdot )$ denote the relative normalized Hilbert polynomial of
a sheaf flat over $S$. 
Let $\Ff _n$ denote the sheaf of $\Lambda$-modules  
$$
\Ff _n := \Ff \otimes _A A/{\bf m}^{n+1}.
$$
It is of pure dimension $d$ on $X$.  Let 
$$
\Ff _n \rightarrow \Gg _n
$$
denote the destabilizing quotient, that is the quotient with the minimal
normalized Hilbert polynomial.  (Note that if $\Gg _0 = \Ff _0$ then $\Ff _0$
is $p$-semistable and we're done---so we assume that this is not the case). We
make the following claims: 
\newline  
1. Let $p_0=p_X(\Gg _0)$.  Then for all $n$, $p_X(\Gg _n)=p_0$.
\newline
2.  There are morphisms $\Gg _n \rightarrow \Gg _{n-1}$ compatible with the
morphisms $\Ff _n \rightarrow \Ff _{n-1}$.
\newline
3.  There there is $q$ such that for $n\geq
q$ we have $\Gg _n \stackrel{\cong}{\rightarrow} \Gg _{q}$.

{\em Proof of 1:} Assume it is known for $n-1$. We have an exact sequence of
$\Lambda$-modules $$
0\rightarrow \Ff _0 \rightarrow \Ff _n \rightarrow \Ff _{n-1} \rightarrow 0.
$$
From this, we obtain a quotient $\Ff _n \rightarrow \Gg _{n-1}$, with
$p_X(\Gg _{n-1})=p_0$.  This shows that the normalized Hilbert
polynomial of the destabilizing quotient $\Gg _n$ is $\leq p_0$. On the other
hand, 
$$
p_X({\rm im}(\Ff _0 \rightarrow \Gg _n)) \leq p_X(\Gg _n)
$$ 
(since $\Gg _n$ is $p$-semistable)---unless this morphism is zero in which
case $\Gg _n = \Gg _{n-1}$ and we're done anyway. 
Therefore (by the definition of $p_0$) we have that $p_0 \leq p_X(\Gg _n)$. 
This proves claim (1).

{\em Proof of 2:} Since $\Gg _{n-1}$ is a quotient of $\Ff _n$ with
$p_X(\Gg _{n-1})=p_0$, it factors through the destabilizing quotient giving the
morphism $\Gg _n \rightarrow \Gg _{n-1}$. 

{\em Proof of 3:} Suppose not. We may assume that $A$ is complete.Let $\Gg :=
\lim _{\leftarrow} \Gg _n$.  This gives  a quotient of $\Ff$ destabilizing
$\Ff$ over the generic point.

Now we proceed with the construction of $\Ff '$. Starting with $\Ff$ (and
assuming that $\Ff _0$ is not $p$-semistable), we construct the 
quotient $\Gg _q$
as above. By statement (3), $\Gg _q$ is the maximal quotient of
$\Ff$ which has normalized Hilbert polynomial $\leq p_0$ and which is supported
over some $Spec (A/{\bf m}^n$. 

Let $\Ff ^{(1)}$ be the kernel of
the map $\Ff \rightarrow \Gg _q$.   Let $p_1$ be the normalized Hilbert
polynomial of the destabilizing quotient $\Qq$ of $\Ff ^{(1)}\otimes _A A/{\bf
m}$. We claim that $p_1 > p_0$.  To see this, suppose to the contrary that
$p_1\leq p_0$. Let $\Kk$ be the kernel of the map  $\Ff ^{(1)}\rightarrow \Qq$,
and let $\Gg '= \Ff /\Kk$.  Then $\Gg '$ is a quotient of $\Ff$ which is an
extension of $\Gg _q$ by $\Qq$; in particular its normalized Hilbert polynomial
is $\leq p_0$.  Furthermore $\Gg '$ is supported over $Spec (A/{\bf
m}^{q+2})$.  This contradicts maximality of $\Gg _q$, showing that $p_1 > p_0$. 

Now start with $\Ff ^{(1)}$ and repeat the same construction to obtain $\Ff
^{(2)}$ etc.; and for each $i$ let $p_i$ be the normalized Hilbert polynomial
of the destabilizing quotient of $\Ff ^{(i)}\otimes _A A/{\bf m}$ (we stop if 
$\Ff ^{(i)}\otimes _A A/{\bf m}$ is $p$-semistable).
We have $p_0 < p_1 < p_2 < \ldots $. Since all of these sheaves are
flat over $S$ (they are subsheaves of $\Ff$ and hence without $A$-torsion), the
Hilbert polynomials of  $\Ff ^{(i)}\otimes _A A/{\bf m}$ are all equal to the
Hilbert polynomial of $\Ff \otimes _AK$.  Finally, the set of $\Lambda$-modules
with a given Hilbert polynomial and with destabilizing quotient having 
normalized
Hilbert polynomial $\geq p_0$, is bounded.  Thus the set of possible
normalized Hilbert polynomials of the destabilizing quotients is finite.  This
shows that the process must
stop.  At the stopping point  $\Ff ^{(i)}\otimes _A A/{\bf m}$ is 
$p$-semistable, and we take $\Ff ' := \Ff ^{(i)}$.
\eop

\subnumero{Application to $M_{Hod}$}

We apply Langton theory to limits of $\Gm$-orbits in $M_{Hod}$. In fact this
applies equally well to the moduli stack $\Mm _{Hod}$.  Suppose $p\in \Mm
_{Hod}(X_s, GL(n))$.  The $\Gm$-orbit of $p$ is a morphism $\Gm \rightarrow
\Mm _{Hod}$ which corresponds to a $\lambda$-connection $(\Ff ', \nabla ')$ on
$X_s\times \Gm $ (where $\lambda : \Gm \rightarrow {\bf A}^1$ is the
projection of the orbit).  

\begin{corollary}
\label{Extn}
With the notations of the above paragraph,  there is an extension
$(\Ff , \nabla )$ of $(\Ff ', \nabla ')$ to a $\lambda$-connection on
$X_s\times {\bf A}^1$, such that $\Ff |_{X_s\times \{ 0\} }$ is a
bundle, is semistable and has vanishing rational Chern classes.
\end{corollary}
{\em Proof:}
First of all note that there exists an extension $(\Ff _1,\nabla _1)$. For this
note that $p$ corresponds to a $\lambda (1)$-connection $(\Ee , \varphi )$ on
$X_s$, and $\Ff '= p_1^{\ast} (\Ee )$ on $X\times \Gm$ with $\nabla = t\varphi
$ (here $t$ denotes the coordinate on $\Gm$).  We can simply put $\Ff
_1=p_1^{\ast}(\Ee )$ and $\nabla _1= t\varphi$ on $X\times {\bf A}^1$.  

Theorem \ref{ThmA} now implies that there exists an extension $(\Ff , \nabla )$
which is semistable over $X_s\times \{ 0\}$.  Note that $\lambda (0)=0$ so the
restriction to $X_s\times \{ 0 \}$ is a Higgs sheaf.  By flatness of $\Ff$ over
${\bf A}^1$, the restriction has vanishing rational Chern classes.  By
(\cite{HBLS} Theorem 2 p. 39), our restriction is actually a bundle. \eop 

\begin{corollary}
\label{limits}
Suppose now that $G$ is any reductive group.
If $p\in M_{Hod} (X/S, G)$ then the limit $\lim _{t\rightarrow 0}t\cdot p
$ exists in $M_{Hod}(X/S, G)$. 
\end{corollary}
{\em Proof:}
We can choose an injection $G\hookrightarrow GL(n)$.  By a variant
of \cite{Moduli} Corollary 9.15 concerning $M_{Hod}$ (we can get this by using
the topological trivialization $M_{Hod} \cong M_{DR} \times {\bf A}^1$ which is
functorial in $G$) the induced map  $M_{Hod} (X_s,G)\rightarrow M_{Hod}(X_s,
GL(n))$ is finite.  Since the limit exists in $M_{Hod}(X_s, GL(n))$ by the
previous corollary, it exists in $M_{Hod}(X_s, G)$. 
\eop

{\bf Question:} What happens for non-reductive groups? If $G={\bf G}_a$ then
the limits again exist (this is exactly the weight property of the Hodge
filtration refered to at the start of the section), so it seems likely that
this will be true in general.

\begin{lemma}
\label{proper}
Suppose $G$ is a reductive group.
Let ${\bf V}\subset M_{Hod}(X, G)$ be the fixed point set of the
$\Gm$-action (note that ${\bf V}$ is concentrated over the origin so in fact
$V\subset M_{Dol}(X, G)$).  Then ${\bf V}$ is proper over $S$.
\end{lemma}
{\em Proof:}
The fixed point set lies over the origin in ${\bf A}^1$, so it is just the 
fixed
point set of the $\Gm$-action on the moduli space of Higgs bundles. For
$G=GL(n)$ this fixed point set is  proper  by (\cite{Hitchin},
\cite{NitsureModuli}, \cite{Moduli} Theorem 6.11). For any $G$, argue
as in the previous corollary.   Alternatively one can obtain properness using
Langton theory as above. 
\eop

\numero{Compactification of $M_{DR}$}

The space $M_{Hod}(X/S,G)\rightarrow {\bf A}^1$ together with the
action of $\Gm$ and the isomorphism between $M_{DR}(X/S,G)$ and the fiber over
$\lambda =1$, allow us to compactify $M_{DR}(X/S,G)$ relative to $S$. This
depends on the properness results of the previous section.

\subnumero{Structure theory for $\Gm$-orbits and construction of some 
quotients}

Suppose $X\rightarrow S$ is a projective morphism with an action of $\Gm$
covering the trivial action on $S$.  Choose a relatively very ample line
bundle $\Ll$ and a compatible action of $\Gm$. Let $V_i$ denote the connected
components of the fixed point set $V$. For each $i$ there is an integer
$\alpha _i$ such that $t\in \Gm$ acts by $t^{\alpha _i}$ on $\Ll |_{V_i}$.

Define a partial ordering $\leq$ on $V$, by saying that $u \leq v$ if there
is a sequence of points $x_1,\ldots , x_m\in X$ with 
$$
\lim _{t\rightarrow 0} x_1 =u
$$
$$
\lim _{t\rightarrow 0} x_k = \lim _{t\rightarrow \infty} x_{k+1}
$$
$$
\lim _{t\rightarrow \infty} x_m =v .
$$
Notice that if $u\leq v$ and $u\in V_i$, $v\in V_j$ then $\alpha _i \geq 
\alpha
_j$ (and if $\alpha _i=\alpha _j$ then $u=v$).

Suppose $V= V_{+} \cup V_{-}$ is a decomposition of
the fixed point set into two disjoint closed subsets (which are consequently
unions of connected components), with the properties that 
$$
v\in V_{+}, \;\; u\in V, \;\; u\leq v \; \Rightarrow \; u\in V_{+}.
$$
and
$$
v\in V_{-}, \;\; u\in V, \;\; u\geq v \; \Rightarrow \; u\in V_{-}.
$$
Put
$$
Y_{+} =  \{ y\in X \; : \;\; \lim _{\lambda \rightarrow \infty}
\lambda \cdot y \in V_{+}\}
$$
and
$$
Y_{-} =  \{ y\in X \; : \;\; \lim _{\lambda \rightarrow 0}
\lambda \cdot y \in V_{-}\} .
$$
These are disjoint closed subsets. They are closed by an argument
similar to the proof of properness in Theorem \ref{ThmD} below. They
are disjoint because if there existed $y\in Y_{+}\cap Y_{-}$ then  
$$
\lim _{\lambda \rightarrow 0}
\lambda \cdot y \; \leq \; 
\lim _{\lambda \rightarrow \infty}
\lambda \cdot y,
$$  
so we would obtain two points $u,v$ with $u\in V_{-}$ and $v\in V_{+}$ but
$u\leq v$; whence $u\in
V_{+}$ and $v\in V_{-}$ (by the conditions on $V_{-}$ and $V_{+}$) 
contradicting the disjointness of $V_+$ and $V_-$. Finally, note that $Y_{+}$
and $Y_{-}$ are, by the nature of their definitions, $\Gm$-invariant.

Let $U:= X-Y_{+}-Y_{-}$. This is a $\Gm$-invariant open set in $X$.  

\begin{theorem}
\label{ThmD}
With the above notations, a universal geometric quotient $U/\Gm$ exists.  It
is separated and proper over $S$.
\end{theorem}

{\em Remark:}  When the subsets $V_+$ and $V_-$ are defined by choosing $a\in
\qq - \zz$ and setting $V_+ = \bigcup _{\alpha _i > a}V_i$ and 
$V_- = \bigcup _{\alpha _i < a}V_i$ then the quotient
defined above is just the geometric invariant theory quotient of the set of
semistable points  (for the linearized action obtained when the 
linearization is
translated by $a$). In particular, in this case the quotient is projective.  I
don't know if the quotient given by this theorem will be projective in
general, nor if it is projective in our example (the compactification of
$M_{DR}$) below.

{\em Proof:}
Let $X^{(pre)}$ denote the set of pre-stable points 
\cite{GIT}. In our case it is easy to see that $X^{(pre)}=X-V$ is just
the complement of the fixed point set. Mumford constructs a universal 
geometric
quotient $\phi :X^{(pre)}\rightarrow X^{(pre)}/\Gm$. This morphism is 
submersive
so $ \phi (U)$ is open, and by the universality we obtain a universal 
geometic
quotient $\phi : U\rightarrow U/\Gm$.  The only problem is to prove that
$U/\Gm$ is separated and proper over $S$.

Suppose $R$ the henselian local ring of $\cc [z]$ at the origin $P$,
with maximal ideal ${\bf m}$ and residue field $R/{\bf m} = \cc$.   
Let $K$ be
the fraction field of $R$, and let $z\in R$ denote a uniformizing parameter.
Let $\tilde{K}$ denote the algebraic closure of $K$ and let $\tilde{R}$ denote
the normalization of $R$ in $\tilde{K}$.  The extension $\tilde{K}$ is obtained
from $K$, as $\tilde{R}$ is obtained from $R$, by adjoining the elements
$z^{1/n}$. Let $\tilde{{\bf m}}$ denote the maximal ideal of $\tilde{R}$.  Note
that $\tilde{R}$ is a valuation ring with $\qq$ as value group, $\tilde{{\bf
m}}$ is the valuation ideal, and $\tilde{R}/\tilde{m} = \cc$.

Any finite extension of $K$ is isomorphic to $K$ (by changing the
parameter). 

Suppose $\eta :Spec (K)\rightarrow U$ is a point.  We have to show that
there is $\varphi \in \Gm (\tilde{K})$ such that $\varphi \eta $ extends to a
point $Spec (\tilde{R})\rightarrow U$, and furthermore that $\varphi$ is
unique up to $\Gm (\tilde{R})$.

The action of $\Gm$ on the point $\eta $ gives a morphism $Spec
(K)\times \Gm\rightarrow X$ which completes to $Spec
(K)\times \pp ^1 \rightarrow X$.  
Let
$$
\eta  _0 := \lim _{t\rightarrow 0} t\cdot \eta  
$$
$$
\eta  _{\infty} := \lim _{t\rightarrow \infty} t\cdot \eta  
$$
as points $Spec (K)\rightarrow X$. These complete to points $Spec
(R)\rightarrow X$. There is a scheme $W$ (the closure of the
graph of the previous morphism in $Spec (R) \times X$) with a diagram  
$$
\begin{array}{ccccc}
Spec (K)\times \pp ^1 & \hookrightarrow & W& \rightarrow &X \\
\downarrow && \downarrow && \downarrow \\
Spec (K)&\hookrightarrow & Spec (R) &\rightarrow &S
\end{array}
$$
where the vertical arrows are proper, and
where $\Gm$ acts compatibly on everything in the top row. The fiber of $W$
over the closed point of $Spec(R)$ decomposes as a string of $\pp ^1$'s
meeting at fixed points for the action.  Let $y_1,\ldots , y_r$ denote the
images in $X$ of the 
fixed points in the string of $\pp ^1$'s over the origin.  Since $W$ was
taken as the closure of the graph, these points are distinct. We can order 
them
so that $y_1= \eta  _0(P)$, $y_r = \eta  _{\infty}(P)$, and $y_i$ is joined to
$y_{i+1}$ by a $\pp^1$ in the fiber.  In this case for a general point $x$ on
the $\pp ^1$ joining $y_i$ to $y_{i+1}$ we have $(\ast )$
$$
\lim _{t\rightarrow 0} t\cdot x = y_i,\;\;\;
\lim _{t\rightarrow \infty} t\cdot x = y_{i+1}.
$$
In particular, refering to our partial ordering above we have
$$
\eta  _0(P) = y_1 < y_2 < \ldots < y_r = \eta  _{\infty}(P).
$$
Note that since
$\eta  \in U(K)$ we have $\eta  _0 \in V_+(K)$ and $\eta  _{\infty} \in V_-(K)$.
Thus $y_1\in V_+$ and $y_r\in V_-$ (as $V_+$ and $V_-$ are closed).  
By the conditions on $V_+$ and $V_-$ there is $k$ such that 
$y_1,\ldots , y_k \in V_+$ and $y_{k+1},\ldots , y_r \in V_-$.
From the definitions of $Y_+$ and $Y_-$ as well as the the property $(\ast )$
we find that the $\pp ^1$ joining $y_i$ to $y_{i+1}$ lies in $Y_+$ if $i<k$ and
in $Y_-$ if $i>k$, whereas it meets $U$ if $i=k$.  The uniqueness of the
$\Gm$-orbit meeting $U$ in the closed fiber gives the separatedness. Choose
$\varphi \in \Gm (\tilde{K})$ so that $\varphi \eta  : Spec (\tilde{K}
)\rightarrow W$ completes to a point $Spec (\tilde{R})\rightarrow W$ with $P$
mapping to a general  point on the $\pp ^1$ joining $y_k$ to $y_{k+1}$.  This
gives the desired $\varphi$ for properness.
\eop

We obtain the following theorem as a corollary.

\begin{theorem}
\label{ThmB}
Suppose $Z\rightarrow S$ is an $S$-scheme on which $\Gm$ acts (acting
trivially on $S$).  Suppose that the fixed point set $W\subset Z$ is proper
over $S$, and that for any $z\in Z$ the limit $\lim _{t\rightarrow 0}t\cdot z
$ exists in $W$. Let $U\subset Z$ be the subset of points $z$ such that the
limit $\lim _{t\rightarrow \infty} t\cdot z$ does not exist in $Z$.  Then $U$
is open and there exists a geometric quotient $Q=U/\Gm$ by the action of
$\Gm$.  This  geometric quotient is separated and proper over $S$.
\end{theorem}
{\em Proof:} Chose a
$\Gm$-linearized very ample line bundle $\Ll$ on $Z$ (this exists by
\cite{GIT}).  Then $\Gm$ acts in a locally finite way on $H^0(Z, \Ll )$ so
we may choose a fixed subspace which gives a projective embedding of $Z$. 
Thus we may assume that $Z\subset \pp ^N$ as a locally closed subscheme, and
that $\Gm$ acts  linearly on $\pp ^N$ preserving $Z$ and inducing the given
action there. Taking the graph we can consider this as an embedding $Z\subset
\pp ^N \times S$.  Let $X$ be the
subscheme closure of $Z$ in $\pp ^N \times S$ (that is, the subscheme
defined by the homogeneous ideal of forms which vanish on $Z$). Note that
$X$ is projective over $S$ and that $\Gm$ acts on $X$ preserving the open set
$Z$. Let $V$ denote the fixed point set in $X$, and let $V_{+}:=W$ be the fixed
point set in $Z$.  Let $V_{-}:= V\cap (X-Z)$ denote the fixed point set in the
complement.  Note that the complement $X-Z$ is closed, hence proper over $S$,
so $V_{-}$ is proper over $S$.  By hypothesis $V_{+}$ is proper over $S$. We
obtain a decomposition $V= V_{+}\cup V_{-}$ as a disjoint union of two closed
subsets. 

Suppose $u,v\in V$ with $v\leq u$.  This means that there is a sequence of
points  $v_0=v , \ldots , v_n = u$ such that $v_i$ is joined to $v_{i+1}$ by a
$\Gm$-orbit (i.e. there is an orbit whose limits are $v_i$ at $\lambda
\rightarrow 0$ and $v_{i+1}$ at $\lambda \rightarrow \infty$).  Suppose $v_i\in
V_{-}$. Then the orbit corresponds to a point $x\in X$ with $\lim _{\lambda
\rightarrow 0}\lambda \cdot x = v_i$.  But if $x\in Z$ then our hypothesis
would give $v_i\in V_{+}$, so $x$ must be in $X-Z$.  Since $X-Z$ 
is closed, the
other limit $v_{i+1}$ must be in $X-Z$ also.  We thus show by 
induction that if
$v=v_0$ is in $V_{-}$ then so is $u=v_n$.  The contrapositive says that if $u$
is in $V_{+}$ then so is $v$.  
We have shown on the one hand that if $v\in V_{-}$ and 
$u\in V$ with $v\leq u$ then $u\in
V_{-}$; and on the other hand that if $u\in V_{+}$ and $v\in V$ with
$v\leq u$ then $v\in V_{+}$.

We are now ready to apply the general construction above.  Define the subsets
$Y_{+}$ and $Y_{-}$ as before, and let $U'$ be the complement $U'=
X-Y_{+}-Y_{-}$. We claim that $Y_{+}$ is the set of points $x\in Z$ such
that $\lim _{\lambda \rightarrow \infty}\lambda \cdot x \in Z$.  
Recall that $Y_{+} := \{ x\in X \; : \;\; \lim _{\lambda \rightarrow \infty}
\lambda \cdot x \in V_{+}\}$.  But if $x\in Z$ with
$\lim _{\lambda \rightarrow \infty}\lambda \cdot x \in Z$ then this limit is in
$V\cap Z = V_{+}$.  On the other hand, if $x\in X$ with  $\lim _{\lambda
\rightarrow \infty} \lambda \cdot y \in V_{+}$ then $x\not \in (X-Z)$ because
$X-Z$ is closed and $V_{+}\cap (X-Z)=\emptyset$. This shows the claim. We next
show that $Y_{-}= X-Z$.  To see this recall that 
$Y_{-} := \{ x\in X \; : \;\; \lim _{\lambda \rightarrow 0}
\lambda \cdot x \in V_{-}\}$. If $x\in Z$ then by hypothesis
$\lim _{\lambda \rightarrow 0}
\lambda \cdot x \in Z$ and $V_{-}\cap Z= \emptyset$, so this shows that
$Y_{-}\subset X-Z$.  On the other hand if $x\in X-Z$ then since $X-Z$ is
closed, $\lim _{\lambda \rightarrow 0}
\lambda \cdot x \in X-Z$ and hence this limit is in $Y_{-}$, which shows that
$Y_{-}=X-Z$.

With the two statements of the previous paragraph we obtain that the
complement $U'$ of $Y_{+}$ and $Y_{-}$ is equal to the subset of points 
of $Z$ whose limits at $\lambda \rightarrow \infty$ do not exist in $Z$, that
is to say that $U'$ is the same as ths subset $U$ described
in the statement of the theorem.  The result of Theorem \ref{ThmD} 
now gives the
universal geometric quotient $U/\Gm$ which is separated and proper over $S$.
\eop

\subnumero{Relative compactification of $M_{DR}(X/S, G)$}

Suppose $G$ is a reductive group and $X\rightarrow S$ a smooth projective
morphism. Apply Theorem \ref{ThmB} to $Z=M_{Hod} (X/S,G)$. The hypotheses 
on $Z$
are given by Corollary \ref{limits} and Lemma \ref{proper}. Note that the open
set $U$ certainly contains the open set   $$
M_{Hod} (X/S,G)\times _{\bf A}^1 \Gm \cong M_{DR}(X/S,G)\times \Gm
$$
as a $\Gm$-invariant open set.  Since the quotient $U\rightarrow Q$ is a
geometric quotient, the image of the open set $M_{DR}(X/S,G)\times \Gm$ is a
geometric quotient of $M_{DR}(X/S,G)\times \Gm$, but we already know the
geometric quotient here, it is just $M_{DR}(X/S,G)$.  Thus our quotient $Q$
contains $M_{DR}(X/S,G)$ as an open set, and $Q$ is proper over $S$.  
We have proved the following theorem.

\begin{theorem}
\label{Th33}
If $G$ is a reductive group and $X\rightarrow S$ a smooth projective morphism,
then there exists a natural relative compactification $\overline{M}_{DR}(X/S,
G)$ proper over $S$ and containing $M_{DR}(X/S,G)$ as an open subset.
\end{theorem}  
{\em Proof:} 
Take $\overline{M}_{DR}(X/S,G)$ to be the quotient
$Q$ of the previous paragraph. \eop

{\em Remark:}
There is
a natural stack-theoretic compactification $\overline{\Mm }_{DR}(X/S, GL(n))$
containing $\Mm _{DR}(X/S,GL(n))$ as an open subset and satisfying the
valuative criterion of properness over $S$. The valuative criterion comes
from Corollary \ref{Extn}. We state this only in the case $G=GL(n)$ because the
finiteness result (\cite{Moduli} Corollary 9.15) used to pass to any group $G$
in the proof of Corollary \ref{limits} is only available for the moduli spaces,
not for the moduli stacks. 

There is a natural Cartier divisor
on $\Mm _{Hod}(X/S,GL(n))$ given as the pullback of the divisor $\{ 0\}
\subset {\bf A}^1$.  This divisor is $\Gm$-invariant, so it projects to a
Cartier divisor in the stack-theoretic quotient $\overline{\Mm
}_{DR}(X/S,GL(n))$; and the open set $\Mm _{DR}(X/S, GL(n))$ is just the
complement of this divisor.  Thus the ``divisor at infinity'' exists as a
natural Cartier divisor.  

In the moduli space compactification $\overline{M}_{DR}(X/S, G)$ the divisor
at infinity is only a Weil divisor.  We can define an intermediate {\em
orbifold compactification} by taking the quotient $M_{Hod}(X/S, G)/\Gm$ in the
sense of stacks.  Here the divisor at infinity is again a Cartier divisor. 
In the orbifold compactification there may be orbifold points corresponding to
fixed points of finite subgroups of $\Gm$.  In the scheme-theoretic
compactification these project to certain quotient singularities.  They
correspond to objects which are like systems of Hodge bundles (or variations
of Hodge structure) except that the Hodge bundles are only indexed by a cyclic
group instead of $\zz$ so the Kodaira-Spencer map $\theta$ can go ``around and
around'' to no longer be nilpotent.

The
orbifold structure of the divisor at infinity is that of the quotient $U\times
_{{\bf A}^1}\{ 0\} /\Gm$.  But this is the quotient of the open subset of
$M_{Dol}(X/S, G)$ corresponding to Higgs bundles with non-nilpotent Higgs
field, by the action of $\Gm$.

\subnumero{Interpretation of the Hodge filtration in terms of the
compactification}

Let $\overline{M}'_{DR}$ denote the orbifold compactification described
above.  There is a principal $\Gm$-bundle over this space, it is just the
total space from before taking the quotient.  This principal bundle corresponds
to a line bundle $\Ll$.  One can see that  $\Ll = 
\Oo _{\overline{M}'_{DR}}(-D)$
where $D$ is the divisor at infinity. Conversely, from the data of
$\overline{M}'_{DR}$ and the divisor at infinity $D$ (which has multiplicity
one) we obtain a line bundle and hence a principal $\Gm$-bundle with a section
defined over $M_{DR}$.  There is only one function on this total space  which
is constant on the multiples of our copy of $M_{DR}$ so this fixes the morphism
to ${\bf A}^1$.  The total space is $M_{Hod}^{\ast}$, the open subset which is
the complement of the locus of Higgs bundles with nilpotent Higgs field.  Thus
we recover most but not all of the Hodge filtration $M_{Hod}$ from our
compactification with its divisor at infinity.

\subnumero{Griffiths transversality revisited}

In the relative case we have obtained a family of compactifications
$$
\overline{M_{DR}}(X/S)\rightarrow S.
$$
On the other hand, recall that $M_{DR}(X/S)$ has the Gauss-Manin connection
which, analytically, translates the fact that $M_{DR}(X/S)^{\rm an}$ is
locally over $S^{\rm an}$ a product of the form $S^{\rm an}\times M_B$ where
$M_B$  is the moduli space of representations of the fundamental group of the
fiber $X_s$.  
The Griffiths transversality condition basically says that the Gauss-Manin
connection has poles of order $1$ along the divisor at infinity.  

In order to make this precise we restrict ourselves to a case where the moduli
space is smooth. Fix a family $X\rightarrow S$ over a base $S$, and suppose $S$
is a smooth curve. Suppose for simplicity that $G$ and Chern class data $c$ are
fixed so that the corresponding unions of components $M_{DR}(X_s, G)_c$ is
smooth and equidimensional (for example $X/S$ a family of curves,
$G=PGL(n)$ and $c$ means we look at bundles of degree $d$ prime to $n$).
We obtain 
$$
M_{Hod} (X/S, G)_c \rightarrow S\times {\bf A}^1.
$$
We claim that this map is smooth.  This is a special case of Conjecture
\ref{triv}, and requires some care. Apply the criterion of (\cite{Hartshorne}
Chapter III Lemma 10.3.A---for which Hartshorne refers to Bourbaki and Altman
and Kleiman) where $t$ is the coordinate on ${\bf A}^1$.  We have to show that
$t$ is not a zero divisor upstairs, and that $M_{Dol}(X/S, G)_c\rightarrow S$ is
flat.  Since all of the fibers of our map are smooth, the only way $t$ could
be a zero divisor is if there were an irreducible component lying over $0\in
{\bf A}^1$.  But Theorem \ref{LocTriv} shows that this is not the case. To
show that $M_{Dol}(X/S, G)_c\rightarrow S$ is flat it suffices (again in 
view of
the smoothness of the fibers) to show that no irreducible
component of $M_{Dol}(X/S,G)_c$  lies over a point in the curve $S$. Any
component of $M_{Dol}(X_s,G)_c$ is contained in the closure of 
$M_{DR}(X_s,G)_c
\times \Gm$, so any component of $M_{Dol}(X/S,G)_c$ is contained in the 
closure
of $M_{DR}(X/S,G)_c \times \Gm$.  If $n$ denotes the dimension of any 
components
of  $M_{DR}(X_s,G)_c$ (by hypothesis these dimensions are all the same) 
then the
dimension of any component of  $M_{Hod}$ is at least $n+2$ and (since $M_{Dol}$
is defined by one equation $t=0$ inside $M_{Hod}$) the dimension of any
component of $M_{Dol}(X/S,G)_c$ must be at least $n+1$.  But the fibers
$M_{Dol}(X_s, G)_c$ are all of dimension $n$, so they cannot contain
irreducible components of $M_{Dol}(X/S,G)_c$. This proves that our map is
flat.  As the fibers are smooth, the map is smooth.

Now we can get back to the thread of our discussion.
Let $U\subset M_{Hod} (X/S, G)_c $ denote
the open set used in defining the compactification.  Let
$\overline{M}'_{DR}(X/S,G)$ denote the orbifold compactification of
$M_{DR}(X/S,G)$ obtained by taking the quotient $U/\Gm$ in the sense of
stacks.  This can introduce orbifold points at places where the stabilizer
is a nontrivial finite subgroup of $\Gm$ (it has to be the $m$-th roots of
unity).  These orbifold points would be replaced by the corresponding cyclic
quotient singularities in the usual compactification defined previously.  The
advantage here is that  $\overline{M}'_{DR}(X/S,G)$ is smooth over $S$.  

Let $D\subset \overline{M}'_{DR}(X/S,G)$ denote the divisor at infinity.  It
is reduced (since $M_{Hod}$ is smooth over ${\bf A}^1$). The Gauss-Manin
connection can be interpreted as a lifting of vector fields on $S$ to vector
fields on $M_{DR}(X/S, G)$. If $p$ denotes the projection to $S$, we obtain a
vector field with coefficients in the line bundle $p^{\ast}(\Omega ^1_S )$
which we denote as
$$
\eta \in H^0(M_{DR}(X/S,G), T(M_{DR}(X/S,G))\otimes p^{\ast}\Omega ^1_S ).
$$ 
Note that it projects to the identity section of $T(S)\otimes \Omega ^1_S=\Oo
_S$. This means that flowing along $\eta$ takes us from one fiber of $p$ to
another.

\begin{theorem}
The Griffiths transversality property says that the vector field $\eta$ has
simple poles along $D$, and the residue is tangent to $D$.  More precisely
let 
$$
\Ff := \frac{T(\overline{M}'_{DR}(X/S,G))\otimes p^{\ast}\Omega ^1_S \otimes
\Oo _D(D) }{T(D) \otimes p^{\ast}\Omega ^1_S \otimes
\Oo _D(D)}
$$
(which is supported on $D$), 
then 
$$
\ker (H^0(\overline{M}'_{DR}(X/S,G), 
T(\overline{M}'_{DR}(X/S,G))\otimes p^{\ast}\Omega ^1_S \otimes
\Oo (D)) \rightarrow H^0(D, \Ff ).
$$
\end{theorem}
{\em Proof:}
We leave this to the reader.
\eop

It should be interesting to study the behavior of the dynamical system given
by this vector field with poles.  The transport between fibers
$M_{DR}(X_s,G)$ and $M_{DR}(X_t,G)$ has the effect of composing the
analytic isomorphisms 
$$
M_{DR}(X_s, G)^{\rm an} \cong M_B(X_s, G)^{\rm an} = M_B(X_t,G)^{\rm an}
\cong 
M_{DR}(X_t, G)^{\rm an}
$$
where the left and right isomorphisms are the Riemann-Hilbert correspondence
and the middle equality comes from the isomorphism of fundamental groups
(which depends on the path we choose from $s$ to $t$).

{\em Exercise:}
Interpret the regularity of the Gauss-Manin connection in way
similar to the above interpretation of Griffiths transversality.
On the moduli space $M(X_{DR/S}(\log D)/S, G)$ the lifts of vector  fields
given by the Gauss-Manin connection will have simple poles along inverse image
of the singular set in $S$.  

\subnumero{Compactifications of spaces of $\Lambda$-modules}

By a technique similar to our construction of the relative compactification of
$M_{DR}$ we have the following general theorem. 

\begin{theorem}
\label{Thm3}
Suppose $X\rightarrow S$ is a projective flat morphism, and $\Lambda$ is a
split almost polynomial sheaf of rings of differential operators on $X/S$. 
Let $M(\Lambda ,P)\rightarrow S$ denote the moduli space of semistable
$\Lambda$-modules with Hilbert polynomial $P$ on $X/S$. Then there exists a
relative compactification, a scheme $\overline{M(\Lambda )}\rightarrow
S$ containing $M(\Lambda )$ as an open set and which is proper over $S$.
\end{theorem}
{\em Proof:}
The proof is the same as the previous one, with the following changes.
We replace the ring $\Lambda _{Hod}$ by the ring $\xi (\Lambda , F)$ on 
$X\times {\bf A}^1$
for the filtration $F$ of $\Lambda$ by degree. Even if not admitted in the
definition of $M(\Lambda ,P)$, we must now admit torsion-free objects in the
space $M(\xi (\Lambda , F), P)$ used to get the compactification.
\eop

We could even obtain the same result with $\Lambda$-modules which are of pure
dimension $d < dim (X/S)$---this is an interpretation of the statement for
$deg (P)=d$.  The proof is again exactly the same.

\subnumero{A total space compactification of $M_{DR}(X/S, GL(n))$} 
Suppose $S$ is smooth, projective, with $S'=S-D$ the complement of a divisor
with normal crossings.  Suppose $X\rightarrow S$ is smooth over $S'$ (we
denote $X':= X\times _SS'$) and has inverse image of $D$ being a divisor with
normal crossings. In this situation
We can get a compactification for the total space $M_{DR}(X'/S', GL(n))$.  
Let $\Lambda = \Lambda _{DR/S (\log D)}$ be the split almost polynomial
sheaf of rings of differential operators corresponding to the formal
groupoid $X_{DR/S}(\log D)$ defined in \S 8.  Apply the construction of
Theorem \ref{Thm3} to obtain a compactification.  We can describe this more
precisely.  There is a formal groupoid $X_{Hod/S}(\log D)$ combining all of
the constructions of \S 8, with underlying scheme $X\times {\bf A}^1$.  Let
$M^{\rm tf} (X_{Hod/S}(\log D), GL(n))$ denote the moduli space for torsion
free semistable objects on $X_{Hod/S}(\log D)$  (with vanishing rational Chern
classes).  Our {\em total space compactification} is
$$
\overline{M}^{\rm tf}(X_{DR/S}(\log D), GL(n)):= U/\Gm 
$$
where $U\subset M^{\rm tf} (X_{Hod/S}(\log D), GL(n))$ is the open set of
points $p$ such that $\lim _{t\rightarrow \infty} tp$ does not exist. It is
proper over $S$ and since $S$ itself is proper, it is compact.

The necessity to
include torsion-free sheaves over the singular fibers is why we must assume
here that the structure group is $GL(n)$. 

Combining the interpretations of Griffiths transversality and regularity of
Gauss-Manin (exercise), we get that the lifts of vector fields on $S$ given
by the Gauss-Manin connection on $M_{DR}(X/S, GL(n))$, have simple poles at all
components of infinity in the total space compactification of $M_{DR}(X/S,
GL(n))$.

\numero{The nonabelian Noether-Lefschetz locus}

A. Beilinson made a comment to the effect that one could get the moduli for
$\zz$-variations of Hodge structure as an intersection between the integral
representations, and the filtered local systems. Of course for $X$ fixed the
moduli space is just a finite set of points, but this becomes interesting when
we let $X$ vary in a family.

For this section we suppose $G=GL(n)$.

Suppose $X\rightarrow S$ is a smooth projective morphism.  Let $V\subset
M_{Dol}(X/S,GL(n))$ denote the fixed point set of the $\Gm$-action; it is also
the moduli space for systems of Hodge bundles or equivalently for complex
variations of Hodge structure \cite{CVHS}.  Let
$V_{DR}$ denote the image in $M_{DR}(X/S,GL(n))$ (note that this is not a 
complex
analytic subset).  On the other hand, let 
$M_B(X_s, GL(n,\zz ))\subset M_B(X_s, GL(n))$ denote the image of $Hom (\pi
_1(X_s,x), GL(n,\zz ))$ (it is the subset of integral representations), and let
$M_{DR}(X/S, GL(n,\zz ))$ denote the
subset of points of $M_{DR}(X/S,GL(n))$ which over each fiber $X_s$ correspond
to elements of $M_B(X_s,GL(n,\zz ))$. Note that $M_{DR}(X/S,GL(n,\zz ))$ 
is a complex
analytic subset of $M_{DR}(X/S,GL(n))$.  Finally put
$$
NL(X/S,GL(n)):= V_{DR}\cap M_{DR}(X/S,GL(n,\zz ) ).
$$
There is a morphism $NL(X/S,GL(n))\rightarrow S$.

\begin{theorem}
For each $s\in S$, the fiber $NL(X/S,GL(n))_s$ is the set of 
isomorphism classes
of integral representations $\rho$ such that $\rho \oplus \rho$ underlies an
integral variation of Hodge structure.  The morphism $NL(X/S,GL(n))
\rightarrow S$ is
proper, and $NL(X/S,GL(n))$ has a unique structure of normal analytic 
variety such
that the inclusions $NL(X/S,GL(n))\rightarrow M_{DR}(X/S,GL(n))$ and
$NL(X/S,GL(n))
\rightarrow M_{Dol}(X/S,GL(n))$ are complex analytic.
\end{theorem}

{\em Remark:}  In this first treatment of the subject, we are ignoring
a possibly more natural non-reduced or non-normal structure of complex
analytic space on $NL(X/S,GL(n))$.

{\em Proof:}  If $\rho \in NL(X/S,GL(n))_s$ then $\rho$ is a complex 
variation of
Hodge structure on $X_s$, and $\rho$ is integral.  It is easy to see that $\rho
\oplus \rho$ has a structure of  integral variation of Hodge structure (see for
example the arguments in  \cite{DeligneMostowVol}, \cite{Hodge1},
\cite{HBLS}). Conversely if $\rho \oplus \rho$ has a structure of variation of
Hodge structure then $\rho$ is fixed by the $\Gm$ action, so $\rho$ lies in 
$V$.

To see that $NL(X/S,GL(n))\rightarrow S$ is proper, it suffices to work locally
near a point $s_0\in S$. By the main result of \cite{DeligneMostowVol}
(done with small variations in the parameters, which still works)
there is a finite subset of points of $M_B(X_{s_0}, GL(2n,\zz ))$ which
correspond to representations are doubles of $n$-dimensional
representations and which could possibly be $\zz$-variations of Hodge structure
on $X_s$ with $s$ in a given relatively compact neighborhood of $s_0$. Over
this neighborhood, $NL(X/S,GL(n))$ is the intersection of $V_{DR}$ (which is
proper over $S$ \cite{Moduli}) with this finite set of sections; thus
$NL(X/S,GL(n))$ is proper over $S$. 

We show the complex analyticity of $NL(X/S,GL(n))$ inside
$M_{DR}(X/S,GL(n))$. The question is local over $S$, so we can fix a
neighborhood $U$ of $s\in S$ and a section $\sigma : U\rightarrow M_{DR}
(X/S,GL(n,\zz )
)$ corresponding to an integral representation of $\pi _1(X_s)$; we have to 
show
that $\sigma ^{-1}(V_{DR})$ is an analytic subset of $U$.
Our representation $\rho$ corresponds to a local system $W$ on $X_s$. 
Let $W_i$ denote the complex irreducible factors of $W$.  Then $\sigma ^{-1}
(V_{DR})$ is the intersection of the subsets $N_i$ of points in $u$ where
$W_i$ admits a complex variation of Hodge structure.  It suffices to show that
$N_i$ are complex analytic.  If $W_i$ does not admit a flat hermitian form
then $N_i$ is empty, so we can assume that it does admit such a form.  This
form $\langle \; ,\; \rangle$ is uniquely determined up to a 
scalar.  Fix an
integer $w$, and data of ranks $r_i$ and degrees $d_i$.  Let 
$$
HF (X/S, W_i, r_i, d_i ) \rightarrow U
$$
denote the parameter scheme for filtrations $F^{\cdot}$ of the local systems
$W_i(s)$ (considered as vector bundles with integrable connection on the
fibers $X_s$) with $r(F^i)=r_i$ and $deg (F^i)=d_i$, and satisfying the
Griffiths transversality condition.  It is an analytic variety over $U$, in 
fact
the pullback of a quasiprojective variety over $M_{DR}(X/S,GL(n))$ by the 
section
$\sigma$ (the parameter variety is the closed subscheme of the Hilbert scheme
of filtrations defined by the Griffiths transversality condition). Let     
$$
HF (X/S, W_i, \langle \; ,\; \rangle , w, r_i, d_i ) \subset 
HF (X/S, W_i, r_i, d_i )
$$
denote the open subset of filtrations such that $F^{\cdot}$
and $F^{\cdot , \perp}$ together determine a Hodge structure of weight $w$ on
the fiber over every $x\in X$. The morphism 
$$
HF (X/S, W_i, \langle \; ,\; \rangle , w, r_i, d_i )\rightarrow U
$$ 
is
injective (because there is at most one structure of complex variation of
Hodge structure on the irreducible representation $W_i$ up to translations,
but the translations are fixed by specifying the ranks $r_i$), so the image is
an analytic subset.   There are only a finite number of possible sets of
degrees
and ranks which can occur (and it suffices to consider $w=0$ for example) so
the union of this finite number of analytic subsets is $N_i$. This shows
that $N_i$ is analytic, hence that the intersection is analytic. Thus 
(taking the union over the finite number of sections we have to consider)    
$NL(X/S,GL(n))$ is an analytic subset of $M_{DR}(X/S,GL(n))$.

From the above construction one gets that over any component of 
$NL(X/S,GL(n))$ the
Hodge filtration of the resulting variation of Hodge structure varies
analytically. Since the associated Higgs bundle is the associated graded of
the Hodge filtration (with $\theta$ as 
the projection of the connection, which
also varies analytically), the associated Higgs bundle varies analytically with
the point in $NL(X/S,GL(n))$, that is to say that $NL(X/S,GL(n))$ is an 
analytic
subset  of
$M_{Dol}(X/S,GL(n))$.
\eop

\begin{corollary}
\label{NLalg}
IF $S$ is projective, then $NL(X/S,GL(n))$ has a structure of normal projective
variety such that the morphisms $NL(X/S,GL(n))\rightarrow M_{DR}(X/S, GL(n))$
and $NL(X/S,GL(n))\rightarrow M_{Dol}(X/S, GL(n))$ are algebraic morphisms.
\end{corollary}
{\em Proof:} We can think of $NL(X/S,GL(n))$ as the normalization of an 
analytic
subvariety of $M_{DR}(X/S,GL(n))$ which is proper over $S$.  In particular 
it is a
closed subvariety of $\overline{M}_{DR}(X/S,GL(n))$ so we can apply GAGA 
to say that
it is algebraic.  Again by GAGA the morphism to $M_{Dol}(X/S,GL(n))$ is 
algebraic.
\eop

After the relatively surprising results of the theorem and this
corollary, the reader might well be asking if $NL(X/S,GL(n))$ isn't just a 
finite
collection of points. In fact, $NL(X/S,GL(n))$ can
have positive dimensional components.  
For example if $Z$ is a surface with a $\zz$-variation of Hodge structure and
then if $X/S$ is a pencil of curves on $Z$, the family of restrictions of the
variation on $Z$ gives a component of $NL(X/S,GL(n))$ dominating $S$. More 
generally
if $X/S$ contains a pencil of curves on $Z$ as a subfamily, then we get a
component of $NL(X/S,GL(n))$ dominating this subfamily.
One might ask whether all positive dimensional components of $NL(X/S,GL(n))$ 
come
from such a construction.

We can also ask for an extension of the result of Corollary \ref{NLalg} to
the quasiprojective case.

\begin{conjecture}
\label{NLalgQPcase}
If $S$ is a quasiprojective variety then $NL(X/S, GL(n))$ is an algebraic
variety and the morphisms to $M_{DR}(X/S, GL(n))$ and 
$M_{Dol}(X/S, GL(n))$ are algebraic.
\end{conjecture}

This conjecture would be the nonabelian analogue of the result of \cite{CDK}. 
It is a globalized version of a problem Deligne posed to me some time ago, of
obtaining a generalization of the finiteness results of
\cite{DeligneMostowVol} uniformly near singularities.  In 
\cite{CDK} it is explained how their result would be a consequence of the
Hodge conjecture.  In a similar way, Conjecture \ref{NLalgQPcase} would be a
consequence of the following nonabelian version of the Hodge conjecture.

\begin{conjecture}
\label{NAHC}
The points of $NL(X_s, GL(n))$ (i.e. the $\zz$-variations of Hodge structure)
are motivic representations on $X_s$.
\end{conjecture}

I don't know who first thought of this conjecture but that must have been a
long time ago.  Of course there is little chance of making any progress on
this---I have presented it only for the light it sheds on Conjecture
\ref{NLalgQPcase}.

The main problem in proving Conjecture \ref{NLalgQPcase} is to make
a local analysis around the singularities of a good completion of the family
$X/S$ to $\overline{X}/\overline{S}$.  We hope that the total
space compactification 
constructed above will be useful for studying this algebraicity question
locally at the singularities of the family.

\numero{An open problem: degeneration of nonabelian Hodge structure}

I would like to end this paper by proposing an open problem for further
research.  This is the problem (motivated at the end of the previous section)
of studying the degeneration of nonabelian Hodge structure which arises from a
degenerating family of varieties $X\rightarrow S$.  It is already complicated
enough to study the case of a degenerating family of curves, so we can suppose
that $X\rightarrow S$ is a family of curves.  Let $0\in S$ denote a point where
the fiber $X_0$ has a simple singularity (node).  Let $U=S-\{ 0\}$ and suppose
$X_U$ is smooth over $U$. We assume $G=GL(n)$ for reasons we will see in a
minute. Then the relative moduli space $M_{DR}(X_U/U , GL(n))$ over $U$ is
provided with the Gauss-Manin connection; with a family of hyperk\"ahler
structures; with a family of Hodge filtrations satisfying Griffiths
transversality; and with a Noether-Lefschetz locus $NL(X_U/U, GL(n))\subset
M_{DR}(X_U/U, GL(n))$ for  which
we would like to prove algebraicity. The problem, then, is to study the
degeneration of all of these structures as one approaches $0\in S$.  

The first step in studying degenerations of variations of Hodge structure in
the abelian case was to have a canonical extension of the underlying
holomorphic bundle, and a regular singularity theorem for the flat connection.
We have obtained the analogues of these things for the nonabelian case.  Note
that in defining a canonical extension, it is important to get something
which is proper over $S$.  This was perhaps not apparent in the abelian case,
where nobody cares about the fact that a vector space is noncompact! But in
the nonabelian case, a noncompact extension would leave room for some
asymptotic behavior as one goes off to infinity, difficult to see.
To get a compactification of $M_{DR}(X_U/U, GL(n))$ we take the total space
compactification $\overline{M}^{\rm tf}(X_{DR/S}(\log 0)/S, GL(n))$ which comes
from looking at the 
moduli of torsion free sheaves with $\lambda$-connection logarithmic at the
singularities, and applying the construction of \S\S 10-11.  After doing 
all of this we obtain a compactification of
$M_{DR}(X_U/U, GL(n))$ which is the analogue of the canonical extension
(together with the Hodge filtration which corresponds to the fiberwise
compactification).  The regularity of the Gauss-Manin connection, coupled with
Griffiths transversality, give that the Gauss-Manin connection on   
$M_{DR}(X_U/U, GL(n))$ over $U$ has poles of order one at infinity (both along
the fiberwise infinity and at the singular point $0\in S$).

We need a good local description of this compactified moduli space near all
points at infinity, specially torsion-free sheaves on the singular fibers but
also at infinity in the fiberwise direction.

One approach to the analysis of the degeneration might be to directly analyze
the lifted vector field giving the Gauss-Manin connection with its simple
poles, and try to deduce asymptotic properties of the transport along this
vector field.

Another option would be to try to analyze the degeneration of the
hyperk\"ahler structure, thinking of it as a family of quaternionic structures
on the fixed variety $M_B(X_s, GL(n))$.  

Among the goals of this study should be: to prove algebraicity of the
Noether-Lefschetz locus, in a nonabelian version of the work of
Cattani-Deligne-Kaplan \cite{CDK}; to obtain a nonabelian version of the
Clemens-Schmid exact sequence; to obtain estimates and asymptotic
expansions for everything in the spirit of the $SL_2$ and nilpotent orbit
theorems; and finally to be able to apply all of this to obtain a devissage
principle for nonabelian cohomology (even with coefficients in higher homotopy
types \cite{kobe}), i.e. a version of the Leray spectral sequence.

\end{document}